\documentclass[a4paper,11pt,twoside]{amsart}
\providecommand{\U}[1]{\protect\rule{.1in}{.1in}}
\usepackage[a4paper,bottom=3.5cm,top=2cm,left=3cm,right=2.5cm]{geometry}

\newtheorem{theo}{Theorem}[section]
\newtheorem{lem}[theo]{Lemma}

\newtheorem{prop}[theo]{Proposition}
\newtheorem{cor}[theo]{Corollary}
\newtheorem{defn}{Definition}[section]

\renewcommand{\H}{{\mathcal H}}
\newcommand{\K}{{\mathcal K}}
\renewcommand{\S}{{\mathcal S}}

\renewcommand{\P}{{\mathbb P}}
\newcommand{\F}{{\mathcal F}}
\newcommand{\G}{{\mathbb G}}
\newcommand{\D}{{\mathcal D}}
\newcommand{\M}{{\mathcal M}}
\newcommand{\B}{{\bf B}}
\newcommand{\E}{{\mathbb E}}
\newcommand{\as}{\mbox{$\P$-a.s.}}
\newcommand{\N}{{\mathbb N}}
\newcommand{\R}{{\mathbb R}}

\newcommand{\X}{{\mathcal X}}

\newcommand{\C}{{\mathcal C}}
\newcommand{\I}{{\bf 1}}
\newcommand{\Z}{{\mathbb {Z}}}
\newcommand{\GG}{{\mathcal G}}

\begin{document}

\title
[Invariance principles under the Maxwell-Woodroofe ]
{ Invariance principles under the Maxwell-Woodroofe condition in Banach spaces}


\author{Christophe Cuny }
\address{Laboratoire MICS, Centralesupelec, Grande Voie des Vignes, 92295 Chatenay-Malabry cedex, FRANCE.}
\email{christophe.cuny@ecp.fr}
\thanks{I am thankful to J\'er\^ome Dedecker for 
providing him with a copy of \cite{Woyczynski}. I wish to thank two anonymous referees for valuable remarks that yielded an improved presentation of the paper}


\begin{abstract}
We prove that, for (adapted) stationary processes, the so-called Maxwell-Wood-roofe 
condition is sufficient for the law of the iterated logarithm   and 
that it is optimal in some sense. That result actually holds in the context of Banach valued stationary processes, including 
the case of $L^p$-valued random variables, with $1\le p<\infty$. 
In this setting we also prove the weak invariance principle, 
hence 
generalizing a result of Peligrad and Utev \cite{PU}. 
The proofs make use of a new maximal inequality and of approximation by 
martingales,  for which some of our results are also new.
\end{abstract}

\maketitle



\section{Introduction} 

Let $(\Omega,\F,\P)$ be a probability space, $\theta$ be an invertible 
bi-measurable measure preserving transformation on $\Omega$ and 
$\F_0\subset \F$ a $\sigma$-algebra such that $\F_0\subset \theta^{-1}(\F_0)$. 
Define a non-decreasing filtration by $\F_n=\theta^{-n}(\F_0)$, for 
every $n\in \Z$ and denote $\E_n:=\E(\cdot |\F_n)$. For every 
$X\in L^1(\Omega,\F,\P)$ write $S_n(X)=X+\ldots +X\circ \theta^{n-1}$.

\medskip

In 2000, Maxwell and Woodroofe \cite{MW} proved the CLT for $(X\circ \theta^n)_{n\ge 0}$ under the condition 
\begin{equation}\label{MWcond}
\sum_{n\ge 1} \frac{\|\E_0(S_n)\|_2}{n^{3/2}}<\infty \,.
\end{equation}
Actually, Maxwell and Woodroofe worked in a Markov chain setting, but in our 
context their condition reads as above.

\medskip

This was a considerable improvement of the martingale-coboundary 
condition of Gordin and Lif\v sic \cite{GL} which in our setting is equivalent to the boundedness of  $(\|\E_0(S_n(X))\|_2)_{n\ge 1}$. 

\smallskip

Moreover, the condition \eqref{MWcond} proved to be useful in applications. It is directly checkable  for linear processes with innovations that are martingale differences, see e.g. Zhao and Woodroofe (Proposition 5 and its proof). 
It leads to the optimal sufficient condition for the CLT in the case of $\rho$-mixing processes, see Merlev\`ede, Peligrad and Utev \cite{MPU} pages 14-15. 
It is implied by the  condition $\sum_n (\log n)^{1+\varepsilon} 
\frac{\|\E_0(S_n)\|_2^2}{n^{2}}<\infty$, which can be checked in the case of 
Markov chains with normal Markov operator, see Cuny \cite{Cuny1}. Finally, 
it is implied by the following condition which is easier to check in 
applications (see e.g. \cite{Cuny} sections 3.1 and 3.3)
\begin{equation}\label{MWstrength}
\sum_{n\ge 1} \frac{\|\E_0(X\circ \theta^{n-1})\|_2}{n^{1/2}} <\infty\, .
\end{equation}
For more situations where the conditions \eqref{MWcond} and \eqref{MWstrength} 
can be checked we refer to \cite{MPU} and the references therein. 

\smallskip

Because of those potential applications several authors tried to have a better understanding of the condition \eqref{MWcond} and its connection with 
probabilistic results such as maximal inequalities, the weak invariance 
principle, the law of the iterated logarithm (LIL) and others. 

\smallskip

A key step toward that better understanding was the paper  
\cite[(2005)]{PU} 
by Peligrad and Utev  who proved 
a new maximal inequality and applied it to deduce the weak invariance principle 
(WIP) under \eqref{MWcond}. Moreover, they proved that \eqref{MWcond} is, in some 
sense, optimal for the CLT.

\smallskip

Later, Peligrad, Utev and Wu  
\cite{PUW} and Wu and Zhao \cite{WZ} proved $L^p$-versions of that 
maximal inequality, in the cases $p\ge 2$ and $1<p\le 2$ respectively and 
obtained new results under $L^p$-versions of \eqref{MWcond}. 

\smallskip

Further extensions of those maximal inequalities have been obtained 
recently by Merlev\`ede and Peligrad \cite{MP}. 

\smallskip

On another hand, the quenched CLT (a strengthening of the CLT), the quenched invariance principle and the law of the iterated logarithm  (LIL)
 have been obtained, 
under various strengthening of \eqref{MWcond}, by Derriennic and Lin 
\cite{DL}, Rassoul-Agha and Sepp\"al\"ainen \cite{RS}, Zhao and Woodroofe \cite{ZW}, Wu and Woodroofe \cite{WW}, Cuny and Lin \cite{CL} 
and Cuny \cite{Cuny1}. 

\smallskip  Very recently, Cuny and Merlev\`ede \cite{CM} investigated 
the martingale approximation method  under $L^p$-versions of 
\eqref{MWcond} and, using a new maximal inequality inspired by \cite{MP}, 
they proved the quenched invariance principle under \eqref{MWcond}. 

\smallskip

In view of all those results, one may expect that \eqref{MWcond} be a 
(sharp) sufficient condition for the LIL,  as well as  for its invariance principle.

\smallskip

In this paper we provide a positive answer to that question 
(the example of Peligrad and Utev \cite{PU} ensures the sharpness). Actually our results
hold in a Banach space setting, including any (separable) $L^p$ spaces 
of a $\sigma$-finite measure space. More precisely, we  prove the almost sure invariance principle 
(ASIP) in $2$-smooth Banach spaces or in $L^p$ spaces with $1\le p<2$. We also 
obtain the WIP for dependent variables taking values in a $2$-smooth Banach space or in a Banach space of cotype 2.

\smallskip

The main motivation for considering Banach-valued variables (especially 
the $L^p$ cases, with $1\le p<\infty$) is the fact that 
there are applications in statistics, in the study of the empirical 
process, see section 6.2. Let us mention some papers in this vein:
 del Barrio, Gin\'e and Matr\'an \cite{BGM}, Berkes, Horv\'ath, Shao and 
 Steinebach \cite{BHSS}, Dedecker and Merlev\`ede \cite{DMASIP} and \cite{DMemp}  or D\'ed\'e \cite{Dede}.  Let us mention also the very recent 
 preprint of Dedecker and Merlev\`ede \cite{DMnew}.

\smallskip

To give a flavour of our results we shall state here a theorem in $L^p$, 
$p\ge 1$. 

\smallskip

Let $(S,\S,\mu)$ be a $\sigma$-finite measure space such that 
$L^1(S,\S,\mu)$ is separable (for instance assume that $\S$ be 
countably generated). Let $X(s)$ be a random variable on 
$(\Omega,\F_0,\P)$  with values in $L^p(S,\S,\mu)$, for 
some $1\le p<\infty$. We shall often consider 
$X$ as a (class of a) measurable function on $(\Omega\times S,\F_0\otimes \S,\P\otimes \mu)$, without mentionning it.

\smallskip

For every integer $n\ge 0$, write $X_n=X\circ \theta ^n$. For every $t\in [0,1]$ and every integer $n\ge 1$, write $S_{n,t}(X):=
\sum_{k=0}^{[nt]-1} X_k+(nt-[nt]) X_{[nt]}$ 
and $T_{n,t}:=S_{n,t}/\sqrt n$.

\medskip

For the sake of clarity, we state the next theorem under a condition in the spirit of \eqref{MWstrength} 
rather than \eqref{MWcond}. With this formulation, the ASIP has already been obtained by the author \cite{Cuny}, when $p=2$; the WIP follows from 
Theorem 3.1 of Dedecker-Merlev\`ede-P\`ene \cite{DMP} (see also Theorem 2.1 of Dedecker-Merlev\`ede-P\`ene \cite{DMPemp}, when $p=2$;
 and the CLT  has been obtained by D\'ed\'e \cite{Dede} when $p=1$. 
 
 
 \medskip
 
 We denote by $\|\cdot \|_2$ the $L^2$-norm on $(\Omega,\P)$.

\begin{theo}
Assume that $\theta$ is ergodic. Let $X \in L^2(\Omega,\F_0,\P,L^p(S))$ ($1\le p<\infty$) be such that $N_{p}(X)<\infty$, where
\begin{gather}
 \label{mcleish0} N_{p} (X)= \sum_{n\ge 1} \frac{\Big(\int_S \big\|\E_0(
 X_{n-1}(s))\big\|_{2}^p \, \mu(ds)\Big)^{1/p}}{n^{1/2}} \qquad \mbox{if 
 $1\le p<2$}\, ,\\
 \label{mcleish}N_{p} (X)= \sum_{n\ge 1} \frac{\big\|\big(\int_S |\E_0(
 X_{n-1}(s))|^p \mu(ds)\big)^{1/p}\big\|_{2}}{n^{1/2}} \qquad \mbox{if 
 $p\ge 2$}\, .
\end{gather}
Then, the process  $((T_{n,t})_{0\le t\le 1})_{n\ge 1}$ converges in law 
in $C([0,1],L^p(S,\mu))$ (to an $L^p(S,\mu$)-valued 
brownian motion); $\big(S_n(X)/\sqrt{nL(L(n))}\big){n\ge 1}$ 
is $\P$-a.s. relatively compact in $L^p(S,\mu)$. Moreover, there exists 
a universal constant $C>0$, such that 
$$
\limsup_{n\to +\infty} \frac{\big(\int_S|\sum_{k=0}^{n-1}X_k(s)|^p\mu(ds)
\big)^{1/p}}{\sqrt{2nL(L(n))}}\le CN_{p}\quad \as
$$ 
\end{theo}
The exact value of the ($\P$-a.s. contant) limsup above may be derived from 
the proof.

\smallskip

Under the assumptions of the theorem an ASIP holds as well, 
see Theorem \ref{theoASIPhilb} and Theorem \ref{theoASIPsmooth}. 

\medskip

Notice that if $(X_n)_{n\ge 0}$ is a  sequence of martingale 
differences (i.e. $\E_{n-1}(X_n)=0$ for every $n\ge 1$) in $L^2(\Omega,L^p(S))$ if $p\ge 2$ 
or in $L^p(S, L^2(\Omega))$ if $1\le p\le 2$ then the condition 
$N_p(X)<\infty$ automatically holds. In this case, 
the WIP  and the ASIP are new when $1\le p< 2$, see section \ref{martsec} 
for  references when $p\ge 2$.

\medskip

When $p> 2$, 
neither the ASIP nor the WIP, can be obtained under condition \eqref{mcleish},  by the method of \cite{Cuny} or \cite{DMP}. Indeed, when $p>$ 2, the only sufficient condition (for the WIP or the ASIP) relying on $(\E_0(X_n))_{n\ge 1}$ that may be derived from the results of \cite{DMP} or \cite{Cuny} is the following (see the proof of Theorem 2.1 of 
\cite{DMPemp} page 758):
$\sum_{n\ge 1} \frac{\big\|\big(\int_S |\E_0(
 X_{n-1}(s))|^p \mu(ds)\big)^{1/p}\big\|_{p}}{n^{1/p}} <\infty$. 

\medskip

Our method of proof follows a classical line. To prove the weak invariance 
principle, we first prove tightness of the underlying process and then 
prove  convergence in law of the finite-dimensional distributions. 
To  prove the almost sure invariance principle (in particular the functional 
law of the iterated logarithm) we first prove a compact law of the iterated logarithm (CLIL) and then invoke an important result of Berger, see Theorem 
\ref{Berger}. The tightness and the CLIL are obtained thanks to suitable maximal inequalities. Our proofs make also use of martingale approximation 
arguments, in particular we first prove all results for 
martingale differences.

\medskip

The paper is organised as follows. In section 2, we recall some definitions 
and lemmas, about probability in Banach spaces, that are necessary for the understanding of the statement and/or the proofs of the results. In section 
3, we state all the results (some of them are new) for martingale with stationnary (and ergodic) increments that are needed in the sequel. 
In section 4, we state maximal inequalities under 
projective conditions. In section 5, we state our limit theorems under projective conditions.  In section 6, we provide several examples including the case of the empirical process. All the results of sections 2-5 are proved in the appendix. 
The fact that our examples satisfy the required conditions is checked in the section 6 itself. 

\medskip

Let us mention that versions of our results may be obtained (with slight modifications)  
for non-adapted stationary processes  or stationary processes 
arising in non invertible dynamical systems.

\section{Generalities on probability on Banach spaces}\label{gen}

Let $(\Omega,\F,\P)$ be a probability space. We will consider Banach-valued random variables. We refer to the book by Diestel and Uhl \cite{DU}  for the basic facts on the topic (definition, conditional expectation...). We shall also use results or notations from Ledoux and Talagrand \cite{LTbook}. In all the paper, we shall be concerned only with separable Banach spaces, in which case  the definitions of a random variable of \cite{DU} and \cite{LTbook} 
co\"\i ncide. Other relevant references on the topic are the books by 
Vakhania, Tarieladze and Chobanyan \cite{VTC} and by Araujo and Gin\'e \cite{AG}.

\medskip

 In all the paper, $(\X,|\cdot |_\X)$ will be a \emph{real} separable  Banach 
 space.  Denote by $L^0(\X)$ the space of (classes modulo $\P$ of) functions from $\Omega$ 
 to $\X$ that are limits $\P$-a.s. of simple (or step) functions.  We define, for every $p\ge1$,  the usual Bochner spaces $L^p$ and their weak versions, as follows

\begin{gather*}
L^p(\Omega,\F,\P,\X)=\{ Z\in L^0(\X)~:~\E(|Z|_\X^p )<\infty \}\, ;\\
L^{p,\infty}(\Omega,\F,\P,\X)=\{ Z\in L^0(\X)~:~ \sup_{t >0} t(\P(|Z|_\X>t))^{1/p} <\infty |\}\, .
 \end{gather*}

For every $Z\in  L^p(\Omega,\F,\P,\X)$, write $\|Z\|_{p,\X}:=(\E(|Z|_\X^p ))^{1/p}$ and for
every $Z\in  L^{p,\infty}(\Omega,\F,\P,\X)$, write $\|Z\|_{p,\infty,\X}:= \sup_{t >0} t(\P(|Z|_\X>t))^{1/p}$.

\medskip

For the sake of clarity, when they are understood, some of the references
to $\Omega$, $\F$ or $\P$ may be omitted. Also, in the case 
when $\X=\R$, we shall simply write $\|\cdot \|_{p}$ or $\|\cdot \|_{p,
\infty}$. 
 Recall that for every $p>1$ there
exists a norm on $L^{p,\infty}(\P,\X)$ (see the proof of 
Lemma \ref{grafakos}), equivalent to the quasi-norm
$\|\cdot \|_{p,\infty,\X}$, that makes $L^{p,\infty}(\P,\X)$ a Banach space.


\medskip


We will state our results in the context of Banach spaces that are $2$-smooth 
or of cotype 2. Let us recall the definitions of those spaces.

\medskip

\begin{defn}\label{smdef}
We say that $\X$ is $2$-\emph{smooth},  if there exists $L\ge 1$,
such that
\begin{equation}\label{smoothdef}
|x+y|_\X^2+|x-y|_\X^2 \le 2(|x|_\X^2+L^2|y|_\X^2) \qquad \forall x,y\in \X\, .
\end{equation}
We shall speak about $(2,L)$-smooth  spaces to emphasize the constant  $L$ 
such that  \eqref{smoothdef}  is satisfied.
\end{defn}
\noindent {\bf Remark.} A Banach space is said to be 
$2$-convex whenever  \eqref{smoothdef} holds in the reverse direction.
\medskip

\begin{defn}
We say that $(d_n)_{1\le n\le N}\subset L^1(\Omega,\F,\P,\X)$ is a 
sequence of  martingale differences, if there exist non-decreasing $\sigma$-algebras $({\mathcal G}_n)_{0\le n\le N}$ such that for every $1\le n\le N$, $d_n$ is 
${\mathcal G}_n$-measurable and $\E(d_n|{\mathcal G}_{n-1})=0$ $\as$ 
\end{defn}

The notion of $2$-smooth Banach spaces is very useful due to the 
inequality \eqref{smooth} below, see for instance Proposition 1 of Assouad \cite{Assouad} (and its corollary).

Assume that  $\X$ is $(2,L)$-smooth. Then,  for every martingale differences $(d_n)_{1\le n\le N}$, we have
\begin{equation}\label{smooth}
\E(|d_1+\cdots +d_N|_\X^2)\le 2L^2 \sum_{n=1}^N\E(|d_n|_\X^2)\, .
\end{equation}


Any Hilbert space is $(2,1)$-smooth. 

\medskip


Any $L^p$ space, $p\ge 2$, (of $\mathbb{R}$-valued functions) associated with a $\sigma$-finite measure is $(2,\sqrt{p-1})$-smooth (see \cite{Pinelis} Proposition 2.1). 




\medskip

We shall also need the concept of Banach spaces of type 2 and of cotype 2. 
These concepts are relevant in the study of the central limit theorem 
in Banach spaces, in particular in their relationship with the notion 
of pregaussian variables that we shall introduce later. 

\medskip

\begin{defn}
We say that a separable Banach space $\X$ is of type $2$ 
(respectively of cotype $2$) if 
there exists $L>0$ such that for every independent random variables 
$d_1,\ldots , d_N\in L^2(\Omega,\X)$, with $\E(d_1)=\ldots=\E(d_N)=0$,
\eqref{smooth}   holds  (respectively, such that 
\eqref{smooth} holds in the reverse direction). 
\end{defn}

Of course, any $2$-smooth Banach space is of type  $2$. 


\medskip

Now, we explain what we mean by an invariance 
principle in a Banach space. 

Let us denote by $\X^*$ the topological dual of $\X$. Let $X \in L^0(\Omega,\X)$ 
be such that for every $x^*\in X^*$, $\E(x^*(X)^2)<\infty$ and $\E(x^*(X))=0$.
We define  a bounded \emph{symmetric} bilinear operator $\K=\K_X$ from $\X^*\times \X^*$
to $\mathbb{R}$, by
\begin{gather*}
\K(x^*,y^*) =\E(x^*(X)y^*(X))  \qquad
\forall x^*,y^*\in \X^* \, .
\end{gather*}
The operator $\K_X$ is called the \emph{covariance operator} associated with $X$.

\begin{defn}\label{gaussiandef}
We say that a random variable $W\in L^0(\Omega, \X)$ is  \emph{gaussian} if, for every  $x^*\in \X^*$,
$x^*(W)$ has a normal distribution. We say that a random variable $X\in L^0(\Omega,\X)$, such that for every $x^*\in X^*$, $\E(x^*(X)^2)<\infty$ and $\E(x^*(X))=0$, 
 is  \emph{pregaussian},
if there exists a gaussian variable $W\in L^0(\Omega,\X)$ with
the same covariance operator,
i.e. such that $\K_X=\K_W$. As in \cite{LTbook}, when $X$ is pregaussian, we shall denote (abusively) by $G(X)$ a gaussian variable having the same covariance operator as $X$.
\end{defn}

\begin{defn}
 We say that a process $(W_t)_{0\le t\le 1}\in 
L^0(\Omega, C([0,1],\X))$ is a Brownian motion with covariance operator $\K$ 
if it is a gaussian process such that for every $x^*,y^*\in \X^*$ and 
every $0\le s,t\le 1$,  
${\rm cov}(x^*(W_s),y^*(W_t))=\min(s,t) \K(x^*,y^*)$.
\end{defn}

\medskip


\begin{defn}
We say that $(X_n)_{n\ge 0}$ satisfies the almost sure invariance principle ({ASIP}) if, without changing its distribution, 
 one can redefine the sequence $(X_n)_{n\ge 0}$ on a new probability space on which  there exists a
sequence  $(W_n)_{n\ge 0}$ of centered i.i.d. gaussian variables, such that
$$|X_0+\cdots +X_{n-1}-(W_0+\cdots +W_{n-1})|_\X=o(\sqrt{nL(L(n))})\qquad \as$$
We say that $(X_n)_{n\ge 0}$ satisfies the weak invariance principle (WIP) 
of covariance operator $\K$
if  $\big( (T_{n,t})_{0\le t\le 1} \big)_{n\ge 1}$ converges weakly 
in $C([0,1],\X)$ to a brownian motion of covariance operator $\K$, where for 
every $t\in [0,1]$ and every $n\ge 1$, $T_{n,t}=S_{n,t}/\sqrt n$ and $S_{n,t}=X_0+\ldots +X_{[nt]-1} +(nt-[nt])X_{[nt]}$.
\end{defn}

\medskip

\begin{defn}
We say that $(X_n)_{n\ge 0}$ satisfies the compact law of the iterated 
logarithm (CLIL) if the sequence $((X_0+\ldots + X_{n-1})/\sqrt{nL(L(n))})_{n\ge 1}$ 
is $\P$-almost surely relatively compact in $\X$. We say that $(X_n)_{n\ge 0}$ satisfies the bounded law of the iterated 
logarithm (BLIL) if the sequence $((X_0+\ldots + X_{n-1})/\sqrt{nL(L(n))})_{n\ge 1}$ 
is $\P$-almost surely bounded in $\X$.
\end{defn}

It has been well known that if $(X_n)_{n\ge 0}$ satisfies the ASIP, it satisfies the CLIL too. However, we have not found a proper reference 
where this is explicitly mentionned, hence we shall provide some arguments. 
Let $(W_n)_{n\in \N}$ be iid gaussian variables taking values in $\X$. 
Then, combining the Theorem page 107 of \cite{Lepage} and Lemma 
3 of \cite{KL} (alternatively, combining Theorem 8.6 and Lemma 3.1 
of \cite{LTbook}), it follows that $(W_n)_{n\in \N}$ satisfies the CLIL. 
Then, the fact that if  $(X_n)_{n\ge 0}$ satisfies the ASIP, it satisfies the CLIL too, readily follows from a standard approximation argument. 

\medskip

It is known (see the discussion page 274 of \cite{LTbook}) that in order to have a central limit theorem 
(or a WIP) for a sequence of  iid $\X$-valued random variables it is 
necessary that the variables be pregaussian.  Hence, to prove invariance principles for stationary sequences, we shall consider only pregaussian 
variables. 

\medskip

\begin{defn}\label{G}
 Let $\G(\Omega,\F,\P,\X)=\G(\X)$  
be the set of pregaussian random variables that are in $L^2(\Omega,\X)$. For 
every $X\in \G(\X)$, denote $\|X\|_{\G(\X)}:=
\|X\|_{2,\X}+\|G(X)\|_{2,\X}$. 
\end{defn}

\begin{lem}\label{lempreg}
Let $\X$ be a real separable Banach space. Then, for every pregaussian variables $X,Y$, the variable $X+Y$ is pregaussian and $\|G(X+Y)\|_{2,\X}\le \|G(X)\|_{2,\X}+
\|G(Y)\|_{2,\X}$. In particular, 
$(\G(\X),\|\cdot\|_{\G(\X)})$, 
is a normed vector space. Actually,  it is a  Banach space. 
\end{lem} 

The proof is given in the appendix. The following result is an obvious 
consequence of Lemma 8.23 of \cite{LTbook}, hence its proof is omitted. 


\begin{lem}\label{LTcot}
Let $\X$ be a real separable Banach space. Let $(\Omega,\F,\P)$ be a probability space, $\theta$ be an 
invertible bi-measurable transformation on $\Omega$. Let $\F'$ be a sub-$\sigma$-algebra of $\F$. Let $X\in L^1(\Omega,\F,\P,\X)$ be pregaussian.
Then $\E(X|\F')$ is pregaussian and for every $n\ge 0$, $X\circ \theta ^n $ is 
 pregaussian. Moreover, $ \|\E(X|\F')\|_{\G(\X)}\le \sqrt 2\|X\|_{\G(\X)}$ and 
$ \|X \circ \theta^n\|_{\G(\X)}\le \sqrt 2\|X\|_{\G(\X)}$.
\end{lem}

\begin{lem}\label{martpreg}
Let $\X$ be real separable Banach space. Let $(\H_n)_{n\ge 1}$ be a non-decreasing filtration and 
let $\H_\infty:= \vee_{n\ge 1}\H_n$. For every $X\in \G(\X)$, 
$\|\E(X|\H_n)-\E(X|\H_\infty)\|_{\G(\X)}\underset{n\to \infty}
\longrightarrow 0$. 
\end{lem}

The proof is given in the appendix.
\medskip

From the above lemmas, we see that it will be very convenient to work 
in $\G(\X)$ in order to obtain invariance principles for a sequence 
$(X\circ \theta^n)_{n\ge 0}$ under conditions involving terms 
of the type $(\E_0(X\circ \theta^n))_{n\ge 0}$.

\medskip

Of course, in 
order to have tractable conditions it is necessary to 
be able to compute $\|X\|_{\G(\X)}$.

\medskip

Let $\X=L^p(S,\mu)$ ($1\le p\le 2$), for some $\sigma$-finite measure, (recall 
that, then, $\X$ is of cotype 2). In this case, the following characterization of
pregaussian variables is part of the folklore. It is due to Vakhania 
 \cite{Vakhania} when $\mu$ is discrete (see \cite[p. 262]{LTbook} for a proof). 
It seems to be essentially due to Rajput \cite{Rajput} for a general $\sigma$-finite measure $\mu$. We provide more details  in the appendix. 
 
\begin{lem}\label{pregLp}
Let $\X=L^p(S,\mu)$ ($1\le p\le 2$), for some $\sigma$-finite measure. 
Then, $X(s)\in L^2(\Omega,\P,\X)$ is pregaussian if and only if 
($X$ is centered and) $\int_S (\E|X(s)|^2 )^{p/2}\mu(ds)<\infty$. 
Moreover, there exists  $C_p>0$, depending only on $p$, such that
\begin{equation}\label{pregLr}
\|G(X)\|_2/C_p\le 
\Big( \int_S (\E|X(s)|^2 )^{p/2}\mu(ds)\Big)^{1/p} \le C_p\|G(X)\|_2\qquad \forall X\in \G(
L^p(\mu))\, .
\end{equation}
Hence, $\G(L^p(\mu))$ may be identified with $\{X\in L^p(S,L^2(\Omega,\R))~:~\E(X)=0\}$.
\end{lem} 
\noindent {\bf Remark.} The above identification makes use of the natural 
embedding of $L^p(S,L^2(\Omega,\R))$ into $L^2(\Omega,\P, L^p(S,\mu))$ (when 
$1<p\le 2$), see 
Lemma \ref{grafakos}.

\medskip
 On another hand, when $\X$ is of type 2, in particular when $\X$ is $2$-smooth, by Proposition 9.24 of \cite{LTbook}, $\|\cdot \|_{\G(\X)}$ is equivalent to $\|\cdot \|_{2,\X}$.

\medskip

Hence, we infer that when $\X=L^p(S,\mu)$, for some $1\le p<\infty$, 
there exists $C_p>0$ such that for every $X\in \G(\X)$,

\begin{gather}
\label{cotype2} \|X\|_{\G(\X)}/C_p\le 
\Big( \int_S (\E|X(s)|^2 )^{p/2}\mu(ds)\Big)^{1/p} \le C_p \|X\|_{\G(\X)} 
\qquad \mbox{if $1\le p\le 2$}\\
\label{type2}\|X\|_{\G(\X)}/C_p\le 
\Big[\E\Big( \int_S |X(s)|^p \mu(ds)\Big)^{2/p}\Big]^{1/2} \le C_p \|X\|_{\G(\X)} 
\qquad \mbox{if $p\ge 2$}
\end{gather}

\bigskip

\medskip

Let us conclude that section with some results concerning the necessity 
of geometric conditions for the WIP, the ASIP or the BLIL, in the case 
of i.i.d. sequences. Those results motivate some of our restrictions in the next sections. 

\medskip

\begin{prop}\label{necessity1}
Let $\X$ be a separable Banach space. Assume that every i.i.d. $\X$-valued  $(X_n)_{n\ge 0}$ in $L^2(\X)$, satisfies the WIP (resp. the ASIP, resp. the 
BLIL). Then, $\X$ is of type 2 (resp. of type 2, resp. of type 
$p$ for every $1\le p<2$). 
\end{prop}

In the case of the WIP, the proposition follows from Theorem 10.5 of \cite{LTbook} (there is even a converse result there). In the case of the BLIL, the result follows from Pisier 
\cite{Pisier1} (see his Remark 2 and the proposition page 208). 
We have no reference for the case of the ASIP, so we provide a proof in the appendix. 

\begin{prop}\label{necessity2}
Let $\X$ be a separable Banach space. Assume that every i.i.d. $\X$-valued 
and pregaussian $(X_n)_{n\ge 0}$ satisfies the WIP (resp. the BLIL). Then, $\X$ is of cotype 2.
\end{prop}

In the case of the WIP, the proposition follows from Theorem 10.7 of \cite{LTbook} (there is even a converse result there). We have no reference for the case of the BLIL, so we provide a proof in the appendix.


\section{The martingale case}\label{martsec}

In this section, we give maximal inequalities and invariance principles
 for martingales with stationary 
differences $(d_n)_{n\ge 0}$. As mentionned in \cite{Cuny}, there is no loss 
of generality in assuming that $d_n=d\circ \theta^n$, where $\theta$ is 
an invertible bi-measurable measure preserving transformation. 
Hence we shall use the notations of the introduction. 

Let us mention that all the results of this section, except the ASIP in 
Proposition  \ref{ASIPmart}, hold for stationary  differences of reverse martingales. Recall that $(d_n)_{n\ge 1}\subset L^1(\Omega,\X)$ 
is a sequence of differences of reverse martingale if 
$\E(d_n|\sigma\{d_k~:~k\ge n+1\})=0$.  

Nevertheless, for stationary sequences of reverse martingales we know that the 
ASIP (as stated in Proposition  \ref{ASIPmart}) holds in the particular case where 
$\X=\R$, see Cuny and Merlev\`ede \cite[Corollary 2.5]{CM}.

 Part of the results stated here are new. We shall discuss their 
 novelty in the sequel.

As mentionned, we use the notations from the introduction.

We first state a maximal inequality that is related to the ASIP.  

\medskip

 For every $X\in L^2(\Omega,\F,\P,\X)$, we consider the
following maximal function
\begin{gather}
\label{maxfun2} \M_{2}(X,\theta,\X):= \sup_{n\ge 1}\frac{|\sum_{k=0}^{n-1}
X\circ \theta^k|_\X}{\sqrt{n L(L(n))}}\, ,
\end{gather}
where $L:=\max(\log,1)$. 

\smallskip

We shall omit the dependence in the parameters $\theta$ and/or $\X$ when they are understood.

\begin{prop}\label{ineLIL}
Let $\X$ be a Banach space. Assume either that $\X$ is a $(2,L)$-smooth 
Banach space or $\X= L^p(S,\S,\mu)$ 
with $\mu$ $\sigma$-finite and $1\le p\le 2$). Then, for every $1<r<2$ 
 there exists $C_r>0$ 
such that 
\begin{equation}\label{mm}
\|\M_2(d)\|_{r,\infty} \le LC_r \|d\|_{\G(\X)}\, .
\end{equation}
\end{prop}
\noindent {\bf Remarks.} Only the case $\X=L^p(S,\S,\mu)$, 
$1\le p<2$ is new here. The proposition is proved in \cite{Cuny} when $\X$ is $(2,L)$-smooth. We do not require $\theta$ to be ergodic. One may wonder 
whether the proposition holds when $\X$ is of cotype 2, or at least $2$-convex, which is an open quetion. 


\medskip



For martingales with stationary and ergodic increments 
in $2$-smooth Banach spaces (admitting a Schauder basis) the CLT has been 
obtained by Woyczy\'nski \cite{Woyczynski}, and the WIP by 
Dedecker-Merlev\`ede-P\`ene \cite{DMPemp} (see the proof of their Proposition 6). Rosi\'nski \cite{Rosinski} 
considered
the case of general arrays of martingale increments a la Brown 
(in the $p$-smooth case). As far as we know, the only CLT for 
martingales  taking values in a Banach space of cotype 2 has been obtained by D\'ed\'e \cite{Dede} in the special case where $\X=L^1(S,\S,\mu)$, 
with $\mu$ $\sigma$-finite. 

Hence, the CLT in the next proposition is only partly new, while the WIP 
seems to be new. Recall that $\G(\F_0,\X)$ has been defined in Definition 
\ref{G}. 

\begin{prop}\label{WIPmart}
Assume that $\theta$ is ergodic. Let $\X$ be a real separable Banach space that is either $2$-smooth 
or of cotype 2. Let $d\in \G(\F_0,\X)$ such that $\E_{-1}(d)=0$. 
Then, $(d\circ \theta^n)_{n\ge 0}$ satisfies the WIP of covariance 
$\K_d$, and there exists $C>0$, such that
\begin{equation}\label{inemart2}
\|\max_{1\le k\le n}|S_k(d)|_\X\|_2\le C n^{1/2}\|d\|_{\G(\X)}
\, .
\end{equation}
\end{prop} 
\noindent {\bf Remark.} The constant $C$ depends only on $\X$.

\begin{prop}\label{ASIPmart}
Assume that $\theta$ is ergodic. Let $\X$ be either a $2$-smooth Banach space or $\X=L^p(S,\S,\mu)$, for 
some $1\le p\le 2$ and $\sigma$-finite $\mu$. Let $d\in \G(\F_0,\X)$ such that 
$\E_{-1}(d)=0$. Then, $(d\circ \theta^n)_{n\ge 0}$ satisfies  the ASIP. Moreover 
\begin{equation}\label{limsup}\limsup_n \frac{|S_n(d)|_\X}{\sqrt{nL(L(n))} }
=
\sup_{x^*\in \X^*,|x^*|_{\X^*}\le 1}\|x^*(d)\|_2 \qquad \as
\end{equation}
\end{prop}
\noindent {\bf Remarks.} {\bf 1.} Since $(d\circ \theta^n)_{n\ge 0}$ satisfies  the ASIP, it satisfies the CLIL too. However, it follows from the proof that the ergodicity of $\theta$ is not necessary for the CLIL. 
As already mentionned the CLIL also holds for stationary differences 
of reverse martingales. \\
{\bf 2.} Only the case $\X=L^p(S,\S,\mu)$, $1\le p<2$ is new here. 
The case where $\X$ is $2$-smooth has been obtained in \cite{Cuny}. 
As in Proposition \ref{ineLIL}, one may wonder whether Proposition \ref{ASIPmart} holds if $\X$ is of cotype 2 or at least 2-convex.


\section{Maximal inequalities under projective conditions}\label{maxsec}

In all of this section we do not require $\theta$ to be ergodic. 

\medskip

Before going further, let us introduce the generalized version 
of the Maxwell-Woodroofe condition that we shall need in the sequel. 
Its relevance will be clear from the next results. 

\medskip

Let $X\in L^2(\Omega,\X)$. Define

\begin{gather}
\label{MW2}\|X\|_{MW_2}:= \sum_{n\ge 0}\frac{\|\E_0(S_{2^n}(X))\|_{\G(\X)}}{2^{n/2}} \,.
\end{gather}
 
To have a better understanding of $\|\cdot \|_{MW_2}$ recall 
that if $\X$ is of type 2 (in particular if $\X$ is $2$-smooth), then 
$\|\cdot\|_{\G(\X)} \le C \|\cdot \|_{2,\X}$ and that if 
$\X=L^r(S,\S,\mu)$ with $1\le r\le 2$ and $\mu$ $\sigma$-finite, 
we have \eqref{cotype2}.

\medskip

In view of applications, let us mention the following easy fact based 
on the observation that $\|\E_0(S_n)\|\le \|\E_0(X)\|+\ldots +\|\E_0(X\circ \theta^{n-1})\| $. There exists $C>0$ such that

\begin{gather*}
\|X\|_{MW_2}\le C \sum_{n\ge 1}\frac{\|\E_0(X\circ\theta^n)\|_{\G(\X)}}{n^{1/2}
} \,.
\end{gather*}

In particular, when $\X=L^p(S,\mu)$, for some $1\le p<\infty$, using \eqref{cotype2} and \eqref{type2}, we see that 
$\|X\|_{MW_2}<\infty$ when $N_p(X)<\infty$, where $N_p(X)$ is defined 
by \eqref{mcleish0} if $1\le p\le 2$ and by \eqref{mcleish} 
if $p\ge 2$. 

We first give an almost sure  maximal inequality, whose proof is based 
on the dyadic chaining in its simplest form, taking into account 
our filtration. Then we derive several other maximal inequalities that 
will be needed later, and that have interest in their own. 

\medskip

There are two important points concerning the following proposition. 
Firstly, it involves the terms $(\E_{-2^k}(S_{2^k}))_{k\ge 0}$ 
which appear in the Maxwell-Woodroofe condition (notice that, by 
Lemma \ref{LTcot}, 
$\|\E_{-2^k}(S_{2^k})\|_{\G(\X)}\le \sqrt 2\|\E_{0}(S_{2^k})
\|_{\G(\X)}$). Secondly, 
for every $k\ge 0$, the sequence  $(d_k\circ \theta^{2^{k+1}\ell})_{\ell\ge 0}$ 
defined below is a stationary sequence of martingale differences. 
The proposition makes use of the following maximal function.  For every $X\in L^1(\Omega,\F,\P,\X)$,  define
\begin{gather}
\label{maxfunp} \M_{1}(X,\theta,\X):= \sup_{n\ge 1}\frac{|\sum_{k=0}^{n-1}
X\circ \theta^k|_\X}{n} \, .
\end{gather}

Recall that, by Hopf's dominated ergodic theorem, see Corollary 2.2 page 
6 of \cite{Krengel}, applied to the real variable $|X|_\X$, we have 
$$
\|\M_1(X,\theta,\X)\|_{1,\infty}\le \|X\|_{1,\X}\, .
$$

\begin{prop}\label{inemax}
Let $X\in L^1(\Omega,\F_0,\P,\X)$. For every $k\ge 0$, write 
$u_k:= |\E_{-2^k}(S_{2^k})|_\X$ and $d_k:= 
\E_{-2^k}(S_{2^k})+\E_{-2^k}(S_{2^k})\circ \theta^{2^k} 
-\E_{-2^{k+1}}(S_{2^{k+1}})$. Then, 
for every integer $d\ge 0$, we have $\P$-almost surely
(with the convention $\sum_{k=0}^{-1}=0$)
\begin{align*}
\max_{1\le i\le 2^d} |S_i|_\X   &  \le \max_{1\le i\le 2^d}\Big|\sum_{\ell=0}^{i-1} (X-\E_{-1}(X))\circ \theta^{\ell } 
\Big|_\X
+\sum_{k=0}^{d-1}\max_{1\le i\le 2^{d-k-1}}\Big|\sum_{\ell=0}^{i-1} 
 d_k\circ \theta^{2^{k+1}\ell}\Big|_\X \\ 
  & 
\qquad \quad +u_d + \sum_{k=0}^{d-1}
\max_{0\le \ell \le 2^{d-1-k}-1}u_k\circ \theta^{2^{k+1}\ell} \, .
\end{align*}
In particular, there exists $C>0$, such that,  
\begin{align}
\nonumber \M_2(X,\theta)  &  \le C \Big( \sum_{k\ge 0} \frac{u_k}{2^{k/2}} 
+ \sum_{k\ge 0} \frac{\big( \M_1(u_k^2, \theta^{2^{k+1}})
\big)^{1/2}}{2^{k/2}} \\
\label{inemax2}  &  \qquad +\M_2(X-\E_{-1}(X), 
\theta) +\sum_{k\ge 0} \frac{\M_2( d_k, \theta^{2^{k+1}}}{2^{k/2}}\Big)\, .
\end{align}
\end{prop}
\noindent {\bf Remark.} That proposition is inspired 
by the works of Peligrad, Utev and Wu \cite{PUW} and of Wu and Zhao  
\cite{WZ}.






\begin{cor}\label{cormax2}
Let $\X$ be Banach space that is either $2$-smooth or of cotype 2. There exists $C>0$ such that for every $X\in \G(\X)$ and every 
integer $d\ge 0$, we have
\begin{gather*}\label{est12}
\|\max_{1\le i\le 2^d} | S_i|_\X\|_2 \le 
C 2^{d/2}\Big( \|X\|_{\G(\X)}+ \sum_{k=0}^d 2^{-k/2}\|\E_{-2^k}(S_{2^k})\|_{\G(\X)}\Big) \, . 
\end{gather*}
In particular, if $\|X\|_{MW_2}<\infty$, then
\begin{gather}\label{inemaxIP}
\sup_{n\ge 1}
\frac{\|\max_{1\le k\le n}|S_k(X)|_\X\|_2}{\sqrt n} 
\le C \|X\|_{MW_2}\,  .
\end{gather}
\end{cor}


\begin{prop}\label{inemaxL2s}
Let $\X$ be either a  $(2,L)$-smooth Banach space or $\X=L^p(S,\S,\mu)$, with 
$1\le p\le 2$. Let 
$X\in L^2(\Omega,\F_0,\P,\X)$ be such that $\|X\|_{MW_2}<\infty$.
For every $1<r<2$, there exists a constant $C_{r}>0$, 
such that
\begin{gather}
\label{inemaxL2} \|\M_2(X)\|_{r,\infty,\X}\le C_{r} \|X\|_{MW_2}\, .
\end{gather}
\end{prop}
\noindent {\bf Remarks.} The constant $C_r$ depends  on $r$ and $L$ if 
$\X$ is $(2,L)$-smooth and on $r$ and $p$ if $\X=L^p(S,\S,\mu)$. 
Define $\|X\|_{H_2}:=\sum_{n\ge 0}\|\E_0(X\circ\theta^n)-\E_{-1}(X\circ \theta^n)\|_{2,\X}
<\infty$. Then, if $\|X\|_{H_2}<\infty$, \eqref{inemaxL2} holds  with 
$\|X\|_{H_2}$ in place of $\|X\|_{MW_2}$. This follows from Theorem 2.10 of \cite{Cuny} when $\X$ is $2$-smooth and,  when $\X=L^r(S,\S,\mu)$,
the proof  may be done exactly  as the proof of Theorem 2.10 of \cite{Cuny}, using \eqref{mm}.

\section{WIP and ASIP under projective conditions} \label{limsec}

In all of this section we DO require $\theta$ to be ergodic. 

\smallskip

We first obtain  martingale approximation results in 
Banach spaces of cotype 2. 

\begin{prop}\label{propWIPcot}
Let $\X$ be a Banach space of cotype 2. Let $X\in \G(\X, \F_0)$ be such that 
$\|X\|_{MW_2}<\infty$. Then there exists $d\in \G(\X, \F_0)$ with 
$\E_{-1}(d)=0$ such that
\begin{equation}\label{wipmax}
\|\max_{1\le k\le n}|S_k(X)-S_k(d)|_\X\|_2 =o(\sqrt n)\, .
\end{equation}
In particular, $(X\circ \theta^n)_{n\ge 0}$ satisfies the WIP of 
covariance operator $\K_d$ and $\K_d(x^*,y^*)=\lim_n 
{\rm cov}(S_n(x^*(X),S_n(y^*(X))/n$ for every $x^*,y^*\in \X^*$.
\end{prop}
\noindent {\bf Remark.} The martingale approximation \eqref{wipmax} 
has been proved in \cite{CM}, 
see Remark 2.4, in the case where $\X$ is a Hilbert space (with an explicit 
expression for $d$). When $\X=\R$ the martingale approximation \eqref{wipmax} is due to Gordin and Peligrad \cite{GP} and the WIP to Peligrad and Utev \cite{PU}.

\begin{theo}\label{theoASIPhilb}
 Let $\X$ be either a Hilbert space or $\X=L^p(S,\S,\mu)$, with 
 $1\le p\le 2$ and $\mu$ $\sigma$-finite. Let  $X\in \G(\X, \F_0)$ be such that 
$\|X\|_{MW_2}<\infty$. Then there exists $d\in \G(\X, \F_0)$ with 
$\E_{-1}(d)=0$ such that
\begin{equation}\label{asipmax}
|S_n(X)-S_n(d)|_\X =o(\sqrt{nL(L(n))}) \qquad \as
\end{equation}
In particular, $(X\circ \theta^n)_{n\ge 0}$ satisfies the ASIP of 
covariance operator $\K_d$ and $\K_d(x^*,y^*)=\lim_n 
{\rm cov}(S_n(x^*(X),S_n(y^*(X))/n$ for every $x^*,y^*\in \X^*$. Moreover 
\begin{equation}\label{kuelbslimsup}
\limsup_n \frac{|S_n|_\X}{\sqrt{2nL(L(n))}}= \sup_{x^*\in \X^*, |x^*|_{\X^*}\le 1} 
\|x^*(d)\|_2\le 10\sqrt 2 \|X\|_{MW_2}\qquad \as
\end{equation}
 \end{theo}
\noindent {\bf Remark.} This result is new even when $\X=\R$. 
In view of the previous proposition, one may wonder whether 
the theorem holds true for Banach spaces of cotype 2 or, at least, for 
$2$-convex Banach spaces.

 \begin{theo}\label{theoASIPsmooth}
Let $\X$ be $2$-smooth Banach space. Let $X\in \G(\X, \F_0)$ be such that 
$\|X\|_{MW_2}<\infty$. Then, $(X\circ \theta^n)_{n\ge 0}$ satisfies the WIP 
and the ASIP of covariance operator $\K$ given by  $\K(x^*,y^*)=\lim_n 
{\rm cov}(S_n(x^*(X),S_n(y^*(X))/n$, for every $x^*,y^*\in \X^*$. 
Moreover, 
\begin{equation}\label{kuelbslimsup2}
\limsup_n \frac{|S_n|_\X}{\sqrt{2nL(L(n))}}= \sup_{x^*\in \X^*, |x^*|_{\X^*}\le 1} 
\big(\K(x^*,x^*)\big)^{1/2}\le 10\sqrt 2 \|X\|_{MW_2}\qquad \as
\end{equation}
\end{theo}

\noindent {\bf Remark.}  Let $\X$ be either as in Theorem \ref{theoASIPhilb} 
or as in Theorem \ref{theoASIPsmooth}. Assume that 
$\|X\|_{H_2}:= \sum_{n\ge 0} \|\E_0(X\circ\theta^n)-\E_{-1}(X\circ \theta^n)
\|_{\G(\X)}<\infty$. Then, $(X\circ \theta^n)_{n\ge 0}$ satisfies the WIP 
and the ASIP of covariance operator $\K$ given by  $\K(x^*,y^*)=\lim_n 
{\rm cov}(S_n(x^*(X),S_n(y^*(X))/n$, for every $x^*,y^*\in \X^*$. 
Moreover, \eqref{kuelbslimsup}  holds with $\|X\|_{H_2}$ in the right-hand side instead of $10\sqrt 2 \|X\|_{MW_2}$. This is proved in  Theorem 2.10 (see also Corollary 2.12) of \cite{Cuny} 
when $\X$ is $2$-smooth and may be proved similarly when $\X =L^r(S,\S,\mu)$ 
using the remark after Proposition \ref{inemaxL2s}.

\medskip

Peligrad and Utev \cite{PU} proved that the condition $\|X\|_{MW_2}<\infty$ 
is optimal (in the sense below) for the CLT. 
Actually,  their example gives also the optimality of the condition 
$\|X\|_{MW_2}<\infty$ for the LIL, see \cite{Cuny-arxiv} for a proof.

\begin{prop}\label{prop-PU}
Let $(a_n)_{n\ge 0}$ be a sequence of positive numbers with $a_n
\to 0$ as $n\to \infty$.  There exist a probability space $(\Omega,\F,\P)$, 
with a transformation $\theta$ and a filtration $(\F_n)_{n\in 
\mathbb Z}$, as in the Introduction, such that there exists 
$X\in L^2(\Omega,\F_0,\P)$ for which 
\begin{equation}\label{MWPU}
\sum_{n\ge 1} a_n \frac{\|\E_0(S_n(X))\|_2}{n^{3/2}}<\infty\, ,
\end{equation}
but $(S_n/\sqrt n)$ is not stochastically bounded and
$$
\limsup_n \frac{|S_n(X)|}{\sqrt{nL(L(n))}} =+ \infty \qquad \as\, .
$$
\end{prop}
\noindent {\bf Remark.} It would be interesting to know whether the condition 
$\sum_{n\ge 1}\frac{\|\E_0(X\circ \theta^n)\|_{2}}{n^{1/2}}<\infty$ is 
also optimal. The optimality of the latter condition for the CLT has been 
recently investigated  by Dedecker \cite{Dedecker}. His arguments do not seem to apply for the LIL.

\section{Examples}

\subsection{A direct example}

We now consider the case of $\rho$-mixing processes for which it is known that the Maxwell-Woodroofe condition is well-adapted, see for instance 
pages 14-15 in  \cite{MPU} or the proof of Lemma 1 (page 548) in 
\cite{PUW}. 

\medskip

Let $(X_n)_{n\in \mathbb{Z}}$ be a stationary $\H$-valued sequence. Define
\begin{gather} \label{defrho}
\rho(n) = \rho ({\mathcal F}_{- \infty}^0, {\mathcal F}_{n}^{\infty}) 
\quad \mbox{and} \quad \psi(n)=\psi ({\mathcal F}_{- \infty}^0, {\mathcal F}_{n}^{\infty})
\end{gather}
where ${\mathcal F}_i^j = \sigma(X_i, \dots, X_j)$ and 
\begin{gather*}
\rho (\mathcal A, \mathcal B) = \sup \Big \{ \frac{{\rm Cov}(X,Y)}{\Vert X \Vert _2 \Vert Y \Vert _2} \, : \, X \in L^2 ({\mathcal A}), Y \in L^2 ({\mathcal B}) \Big \} \, ;\\
\psi (\mathcal A, \mathcal B) = \sup \Big \{\frac{|\P(A\cap B)-\P(A)\P(B)|}{
\P(A)\P(B)} \, : \, A\in {\mathcal A}, B\in({\mathcal B}) \Big \} \,
\end{gather*}

It is well-known that $\rho(n)\le \psi(n)$, see for instance Proposition 3.11 page 76 of \cite{Bradley}.

\medskip

We have 

\begin{cor}\label{rhomixing}
 Assume that
\begin{gather}
\label{mixcond}\sum_{n\ge 1} \rho(2^n) <\infty 
\end{gather}
Then, $\|X\|_{MW_2}<\infty$. 
\end{cor}
\noindent {\bf Remarks } The  condition $\rho(2^n)=
O(1/n^{1+\varepsilon})$ has been proven to be sufficient  in \cite{Shao} 
(for any $\varepsilon>0$), when $\H=\mathbb{R}$. The sufficiency of 
\eqref{mixcond} has been obtained very recently by Lin and Zhao 
\cite{LZ}, when $\H=\R$. \\
 Sharipov \cite{Sharipov} obtained the conclusion of 
the corollary under the condition $\sum_n \psi(n) <\infty$. However he 
assumes weaker moment conditions and the variables are allowed to 
take values in a $2$-smooth Banach space.\\ 
\noindent {\bf Proof. } It suffices to prove that
\begin{equation}\label{dy}
\sum_n \frac{\|\E_0(S_{2^n}(X_0))\|_{2,\H}}{2^{n/2}}<\infty \, ,
\end{equation}

Let $(e_i)_{i\ge 0}$ be an orthonormal basis of $\H$, and write 
$Y_0^{(i)}:= \langle X_0, e_i\rangle_\H$. We have 
$$
\|\E_0(S_{2^n}(X_0))\|_{2,\H}^2 =\sum_{i\ge 0} \E\big[ 
\big(\E_0(S_{2^n}(Y_0^{(i)} ))\big)^2\big]\, .
$$
 Now, 
it follows from the computations page 15 of \cite{MPU} combined with 
Lemma 3.4 of \cite{Peligrad} that 
$$
\E\big[ 
\big(\E_0(S_{2^n}(Y_0^{(i)} ))\big)^2\big] \le C \E((Y_0^{(i)})^2)
\big(\sum_{k=0}^n 2^{k/2}\rho(2^k)\big)^2\, .
$$
Using that $\sum_{i\ge 0}(Y_0^{(i)})^2= |X_0|_\H^2$ 
we see that \eqref{dy} is satisfied as soon as 
$$
\sum_n \frac{1}{2^{n/2}} \sum_{k=0}^n 2^{k/2} 
\rho(2^k) <\infty\, ,
$$
which holds, by \eqref{mixcond}.   \hfill $\square$

\subsection{Applications to the empirical process}

Let $(\Omega,\F,\P)$ be a probability space, $\theta$ be an invertible 
bi-measurable measure preserving transformation on $\Omega$ and 
$\F_0\subset \F$ a $\sigma$-algebra such that $\F_0\subset \theta^{-1}(\F_0)$. 
Define a non-decreasing filtration by $\F_n=\theta^{-n}(\F_0)$, for 
every $n\in \Z$ and denote $\E_n:=\E(\cdot |\F_n)$. 

\medskip

Let $Y\in L^0(\Omega,\F_0,\P)$. For every $n\in \mathbb{Z}$, let $Y_n:= Y\circ \theta^n$ and $X_n:= t\mapsto \I_{Y_n\le t}-F(t)$,
where $F(t)=\P(Y\le t)$.

\medskip

Let $p\ge 1$. For every $\sigma$-finite Borel measure $\mu$ on 
$\mathbb{R}$, we may see
$(X_n)_{n\in \mathbb{Z}}$ as a process with values in the  Banach space $L^p(\mathbb{R},\mu)$ (which is $2$-smooth when $r\ge 2$), as soon as
\begin{equation}\label{condminCM}
\int_0^\infty (1-F(t))^p \mu(dt)\quad +\quad \int_{-\infty}^0 F(t)^p\mu(dt)
<\infty \,,
\end{equation}
which is satisfied whenever $\mu$ is finite.

\medskip

Define $F_\mu$ by $F_\mu(x)=-\mu([x,0[)$ if $x\le 0$ and $F_\mu(x)=
\mu([0,x[)$ if $x\ge 0$.  Then, under \eqref{condminCM},
$X_0\in L^2(\Omega,L^p(\mu))$
if and only if
\begin{equation}\label{condminCMbis}
\E(|F_\mu(Y_0)|^{2/p})<\infty \, .
\end{equation}

We want to understand the asymptotic behaviour of the process $F_n=S_n(X)/n$ (with values in $L^2(\Omega,\F_0,\P,L^p(\mathbb{R},\mu))$), and more particularly of $D_{n,p}(\mu):=\|F_n\|_{p,\mu}$. 

\medskip

Notice that when $\mu$ is the Lebesgue measure $\lambda$ 
and $p=1$, $D_{n,1}(\lambda)$ represents the Wasserstein distance 
between the empirical distribution and the true distribution.


\medskip

Let us introduce some dependence coefficients. For every 
$Y\in L^1(\Omega,\F,\P)$ and every $1\le p\le \infty$, 
define 

\begin{gather*}
\check \tau_{\mu,p}(\F_0,Y_n):= \Big \| \Big(\int_\R \big|\, 
\P(Y_n\le t|\F_0)-F(t)\, \big|^p 
\mu(dt)\big)^{1/p} \, \Big\|_2\qquad \mbox{if $p\ge 2$}\, ,\\ 
\check \tau_{\mu,r}(\F_0,Y_n):=  \Big(  \int_\R \big\|\, 
\P(Y_n\le t|\F_0)-F(t)\, \big\|_2^p 
\, \mu(dt)\Big)^{1/p}  \qquad \mbox{if $1\le p <2$}\, .
\end{gather*}

 When $p\ge 2$, $\check\tau_{\mu,p}(\F_0,Y_n)=\tau_{\mu,p}(\F_0,Y_n)$, 
where $\tau_{\mu,p}(\F_0,Y_n)$ appears for instance in \cite{DMemp} 
(notice that our notations are slightly different). 

\medskip


Let us notice that both \eqref{condminCM} and \eqref{condminCMbis} are satisfied 
as soon as $\tau_{\mu,p}(\F_0,Y_0)<\infty$.

\medskip

\begin{theo}\label{emp}
Let $Y\in L^0(\Omega,\F_0,\P)$ and $(S,\S,\mu)$ be a $\sigma$-finite 
measure space. Let  $1\le p<\infty$. Assume that 
$$
\sum_{n\ge 0} \frac{\check \tau_{\mu,p}(\F_0,Y_n)}{n^{1/2}}<\infty\, .
$$
Then $(X_n)_{n\ge 1}$ satisfies the WIP and 
the ASIP. In particular, $(n^{1/2} D_{n,p})$ converges in law to an 
$L^p$-valued gaussian variable, with covariance operator given by 
$\K_\mu(f,g)$ and 
$$
\limsup_n \frac{n^{1/2}}{\sqrt{2L(L(n))}}D_{n,p}(\mu) = 
\Lambda_\mu \qquad \as \, ,
$$
for some $\Lambda_\mu\ge 0$. \\
Let $p'$ be the conjugate of $p$. 
 We have $$\K_{\mu}(f,g)=\lim_{n\to +\infty} \E\big( \int_S f(s)S_n(s) \mu(ds) \int_S g(t)S_n(t) 
\mu(dt)\big)/n\qquad \forall f,g \in L^{p'}(S)\, ,$$ and $\Lambda_{\mu,p}=
\displaystyle \sup_{\|f\|_{p',\mu}\le 1}\Gamma_{\mu,p}(f)$ where 
$\Gamma_{\mu,p}(f)=\lim_n \|\int_S f(s)S_n(s)\, \mu(ds)\|_2/\sqrt n$. 

\end{theo}
\noindent {\bf Remark.} Actually, if $p'$ denotes the conjugate of $p$, 
 we have $\Lambda_{\mu,p}=
\displaystyle \sup_{\|f\|_{p',\mu}\le 1}\Gamma_{\mu,p}(f)$ where 
$\Gamma_{\mu,p}(f)=\lim_n \|\int_S f(s)S_n(s)\, \mu(ds)\|_2/\sqrt n$. 
Since Theorem \ref{emp} is a straightforward application of the results of 
section \ref{limsec}, we omit the proof.

\medskip

In a series of paper, Dedecker and Merlev\`ede obtained the WIP or the ASIP under conditions on the coefficients 
$\check \tau_{\mu,p}$, when $p\ge 2$. In  \cite{DMemp} they  
studied the WIP and in \cite{DMASIP} the ASIP. When $p>2$, 
 their results rely on a  condition a la Gordin, hence yield to stronger conditions than ours. When $p=2$, they use a very different approach and their results have different range of applicability. 
 
When $p=1$,  D\'ed\'e \cite{Dede} obtained the CLT under the same condition as 
above. 

\medskip

In order to apply Theorem \ref{emp} we shall further study the coefficients $\check 
\tau$, and estimate them thanks to other coefficients that are known to be  computable in many situations (see e.g. Dedecker and Prieur \cite{DP}). 

\medskip

Let us define the coefficients $\tilde \phi$ and $\tilde \alpha$, as 
defined in Dedecker and Prieur \cite{DP}. For 
every $n\ge 1$, define  
\begin{gather*}
\tilde \phi(n):=\sup_{t\in \R} \|\P(Y_n\le t|\F_0)-F(t)\|_\infty\, \\
\tilde \alpha (n):= \sup_{t\in \R}\|\P(Y_n\le t|\F_0)-F(t)\|_1\, .
\end{gather*}

\begin{lem}\label{lemmix1}
Assume that $\mu$ is finite. Let $p\ge 1$ and define $q:=\max(2,p)$. 
For   every $n\ge 1$, we have
\begin{gather*}
\tau_{\mu,p}(\F_0,Y_n)\le \mu(\R)^{1/p}\tilde \phi(n) \, ,\\
\tau_{\mu,p}(\F_0,Y_n)  \le \mu(\R)^{1/p} \tilde \alpha (n)^{1/q}\, .
\end{gather*}
\end{lem}
\noindent {\bf Proof.} The first inequality is obvious. The second one follows 
from the fact that for every $s\ge 1$, $\big\|\, \P(Y_n\le t|\F_0)-F(t)\, 
\big\|_s\le \big\|\, \P(Y_n\le t|\F_0)-F(t)\, \big\|_1^{1/s}$. \hfill $\square$


\medskip

\begin{lem}\label{lemmix2}
 Let $1\le p\le 2$.
For   every $n\ge 1$, we have
\begin{gather}
\label{phi}\check\tau_{\mu,p}(\F_0,Y_n)\le \sqrt 2 \Big(\int_0^\infty 
 \big(F(t)(1-F(t))\big))^{p/2}\mu(dt)\Big)^{1/p}\tilde \phi(n)^{1/2}  \, .\\
\label{alpha}\check\tau_{\mu,p}(\F_0,Y_n)  \le \sqrt 2 \Big(\int_0^{+\infty} \, 
\Big( \min \big[\tilde \alpha_n, F(t)(1-F(t)\big]\Big)^{p/2}\, 
\mu(dt)\Big)^{1/p}\, .
\end{gather}
\end{lem}
\noindent {\bf Proof.} Notice  that, for every $t\in \R$,
\begin{gather*}
\label{mixcoeff3} \big\|\, 
\P(Y_n\le t|\F_0)-F(t)\, \big\|_2 ^2 \le 
2 \tilde \phi(n)(1-F(t))F(t)\
\, .
\end{gather*}
Hence, $\eqref{phi}$ follows.

\medskip
 Using that for every $t\in \R$, 
\begin{align*} 
  &  \big\|\, \P(Y_n\le t|\F_0)-F(t)\, \big\|_2 ^2
\le \big\|\, \P(Y_n\le t|\F_0)-F(t)\, \big\|_1\le \tilde 
\alpha(n)\, ,\\
\mbox{and } \quad  &  \big\|\, \P(Y_n\le t|\F_0)-F(t)\, \big\|_2^2
\le 2F(t)(1-F(t))\, ,
\end{align*}
we see that \eqref{alpha} holds. \hfill $\square$




\begin{theo}\label{theomix}
Let $1 \le p\le 2$. Assume either of the following items. 
\begin{itemize}
\item [$(i)$] $\int_0^\infty 
 (F(t)(1-F(t)^{p/2}\mu(dt)\, <\infty$ and $\sum_{n\ge 1} 
n^{-1/2}\tilde \phi(n)^{1/2} <\infty$.
\item [$(ii)$]  $\mu=\lambda$  the Lebesgue measure and  $\sum_{n\ge 1} n^{-1/2}\big(\int_0^{\tilde\alpha(n)}x^{p/2-1}Q(x) \, dx\big)^{1/p}<\infty$, where $Q(x):=\inf\{t\ge 0~:~\P(|Y|>t)\le x\}$.
\end{itemize}
Then, the conclusion of Theorem \ref{emp} holds.
\end{theo}
\noindent {\bf Remark.} A better sufficient condition, in terms of $(\tilde \alpha(n))$ for the WIP has been obtained by 
Dedecker and Merlev\`ede \cite{DMnew} when $p=1$, see their sections 4.4 and 5. 

\smallskip

\subsection{Proof of Theorem \ref{theomix}} 
The conclusion under $(i)$ follows from Theorem \ref{emp} and Lemma 
\ref{lemmix2}. To prove item $(ii)$, in view of Theorem \ref{emp} and Lemma 
\ref{lemmix2}, it suffices to prove that (notice that $F(t)(1-F(t))\le \P(|Y|\ge |t|)=
\P(|Y|>|t|)$ for $\lambda$-a.e. $t\in \R$)
$$
\sum_{n\ge 1} \frac1{n^{1/2}}\Big(\int_0^{+\infty} \, \Big( \min \big[\tilde \alpha_n, (\P(|Y|>t))\big]\Big)^{p/2}\, 
\lambda(dt)\Big)^{1/p} <\infty\, .
$$
Now, 
\begin{equation}\label{jjj}
\int_0^{+\infty} \, \Big( \min \big[\tilde \alpha_n, (\P(|Y|>t))\big]\, dt
\le \tilde \alpha (n)^{p/2} Q(\tilde \alpha(n))+  \int_{Q(\tilde 
 \alpha(n))}^{+\infty} (\P(|Y|>t))^{p/2-1} \, dt
 \end{equation}
Since $Q$ is non-increasing, we see that $(ii)$ implies that 
$$
\sum_{n\ge 1}\frac{\tilde \alpha (n)^{1/2} \big(Q(\tilde \alpha(n))\big)^{1/p}}
{n^{1/2}}<\infty\, ,
$$ 
hence, it remains to deal with the second term in the right-hand side of \eqref{jjj}. 

\medskip

We have 
\begin{gather*}
\int_{Q(\tilde 
 \alpha(n))}^{+\infty} (\P(|Y|>t))^{p/2-1} dt=\int_{Q(\tilde 
 \alpha(n))}^{+\infty} \Big( \int_0^1 \frac{p}2x^{p/2-1} {\bf 1}_{\{x\le \P(|Y|>t)\}}
 dx \Big)\, dt\\
 \le \int_0^{\tilde \alpha (n)} \frac{p}2x^{p/2-1}\Big(\int_0^{Q(x)}dt \Big) \, 
 dx =\int_0^{\tilde \alpha (n)} \frac{p}2x^{p/2-1}Q(x)\, dx\, ,
\end{gather*} 
 and the proof is complete. \hfill $\square$

 \begin{appendix}

 \section{Proof of the results of section \ref{gen}}

\subsection{Proof of Lemma \ref{lempreg}} 

\medskip

Let $X,Y \in \G(\X)$. Consider the Banach space $\C:=\X\times \X$ with norm 
$|(x,y)|_\C  :=(|x|_\X^2+|y|_\X^2)^{1/2}$. Let us prove that 
$(X,Y)\in \G(\C)$. Let $G(X)$ and $G(Y)$ be independent gaussian variables 
with same covariance operator as $X$ and $Y$ respectively. Then, 
$(G(X),G(Y))$ is a gaussian variable taking values in $\C$. Now, 
for every $x^*,y^* \in \X^*$, we have
$$
\E((x^*(X)+y^*(Y))^2)\le 2\E\big[ \big(x^*(G(X))\big)^2+\big(y^*(G(Y))\big)^2
\big]=2 \E((x^*(G(X))+y^*(G(Y)))^2)\, .
$$
Hence, by Lemma 9.23 of \cite{LTbook}, $(X,Y)\in \G(\C).$ Let $(U,V)$ be a gaussian 
variable with values in $\C$ with same covariance operator 
as $(X,Y)$. Clearly, $U+V$ is gaussian and has same covariance operator 
as $X+Y$. Hence, $X+Y$ is pregaussian and we may take $G(X+Y)=U+V$. 
Similarly, we may take $G(X)=U$ and $G(Y)=V$. Now, 
\begin{gather*}
\|G(X+Y)\|_{2,\X}=\|U+V\|_{2,\X}\\
\le \|U\|_{2,\X}+\|V\|_{2,\X}=
\|X\|_{2,\X}+\|Y\|_{2,\X}\, .
\end{gather*}
Hence, $\|\cdot \|_{\G(\X)}$ is a norm on $\G(\X)$.
\medskip

Let us prove that $\G(\X)$ is a Banach space.
\smallskip

Let $(X_n)_{n\ge 1}$ be Cauchy in 
$(\G(\X)),\|\cdot \|_{\G(\X)})$. Hence, 
$(X_n)_{n\ge 1}$ is Cauchy in $L^2(\Omega,\X)$, so it  converges, say to $X$ in 
$L^2(\Omega,\X)$. We just have to prove that $X$ is pregaussian and 
that $(X_n)_{n\ge 1}$ admits a subsequence converging to $X$ for 
$\|\cdot \|_{\G(\X)}$. By assumption, there exists a subsequence 
$(X_{n_k})_{k\ge 1}$ such that $\|X_{n_k}-X_{n_{k+1}}\|_{\G(\X)}\le 2^{-k}$. 
Then $X=-X_{n_1}+\sum_{k\ge 1}X_{n_k}-X_{n_{k+1}}$ with convergence in $L^2(\Omega, \X)$.

\medskip

Extending our probability space, if necessary, we may assume that there exists 
 a sequence $(G_k)_{k\ge 0}$ of independent gaussian variables taking values in $\X$, such that $G_0=G(X_{n_1})$ and for every $k\ge 1$, 
 $G_k= G(X_{n_{k+1}}-X_{n_k})$. Then, $G:=\sum_{k\ge 0} 2^{k/2}G_k$ defines a 
 gaussian variable. Moreover, for every $x^*\in \X^*$, we have, using 
 Cauchy-Schwarz, 
 \begin{gather*}
 \E(x^*(X)^2)= \E\big[ \big(x^*(-X_{n_1})+\sum_{k\ge 1} x^*(X_{n_{k+1}}-X_{n_k}
 )\big)^2\big] 
 \\\le 2 \Big( \E((x^*(-X_{n_1}))^2)+ \sum_{k\ge 1} 2^k\E 
 ((x^*(X_{n_{k+1}}-X_{n_k})^2) \Big) 
 = 2 \E\big[ \big(x^*(\sum_{k\ge 0} 
 2^{k/2}G_k)\big)^2\big]\, .
 \end{gather*}
It follows from Lemma 9.23 of \cite{LTbook} that $X$ is pregaussian. By a similar argument, 
using the second half of Lemma 9.23 of \cite{LTbook},  we see that $\E(|G(X-X_{n_m})|^2
\to 0$ as $m\to +\infty$, and the proof is finished. \hfill $\square$

\subsection{Proof of Lemma \ref{martpreg}}  Let $x^*\in \X^*$. Clearly,  we may assume that $X$ is $\H_\infty$-measurable. Denote $X_n:=\E(X|\H_n)$. Then $(X_n)_{n\ge 1}$ is a martingale 
converging in $L^2(\Omega,\X)$ to $X$ (see for instance Proposition V.2.6. 
of Neveu \cite{Neveu}. It suffices to 
prove that $\|G(X-X_n)\|_{2,\X}$ converges to $0$. Using Lemma \ref{LTcot},  we have

\begin{gather*}
\E[( x^*(X_n-X))^2]\le 2(\E[( x^*(X))^2]+\E[( x^*(X_n))^2])  \le 
6 \E[( x^*(G(X)))^2]\, .
\end{gather*}

\medskip

Since $X_n-X$ is (clearly) pregaussian,  we infer that 
$$
\E[( x^*(G(X_n-X)))^2]\le  6 \E[( x^*(G(X)))^2]\, .
$$
Then, it follows from the discussion pages 73-74 of \cite{LTbook}, that 
$(G(X_n-X))_{n\ge 1}$ is tight, hence converges in probability to 0, 
since for every $x^*\in \X^*$,  $(x^*(G(X_n-X)))_{n\ge 1}$ 
converges in probability to 0 (recall that $\|x^*(G(X_n-X))\|_2
=\|x^*(X_n-X)\|_2\underset{n\to\infty}\longrightarrow 0$. 

\medskip

Let $\varepsilon >0$. There exists $n_\varepsilon\ge 1$ such that 
$\P(|G(X_{n_\varepsilon}-X)|_\X >\varepsilon)<1/2)$.  
In particular, the median of the gaussian variable $G(\tilde 
X_{n_\varepsilon}-X)$ 
is smaller than $\varepsilon$, and it follows from the last assertion of 
Lemma 3.2 of \cite{LTbook}, that there exists a universal $C>0$ such that 
$\|G(\tilde X_{n_\varepsilon}-X)\|_{2,\X}\le C\varepsilon ^2$, 
and the proof is finished. \hfill $\square$

\subsection{Proof of Lemma \ref{pregLp}}  

Let $X(s)\in L^2(\Omega,\P,L^p(S,\mu))$ be pregaussian. Hence there exists  
a gaussian variable $W$ on $(\Omega,\F,\P)$ with values in $L^p(S,\mu)$ with same 
covariance operator than $X$. By Theorem 3.1 of 
Rajput \cite{Rajput}, we may see $W$ as a gaussian process 
$(W(s))_{s\in S}$ whose paths are $\P$-a.s. in $L^p(S,\mu)$. Then, 
\begin{gather*}
\infty > \|G(X)\|_{2,L^p(\mu)} =\|W\|_{2,L^p(\mu)} \ge C_p \|W\|_{p,L^p(\mu)}
=C_p\Big( \int_S \E(|W(s)|^p)\, \mu(ds)\Big)^{1/p} \\
=\tilde C_p \Big( \int_S (\E(|W(s)|^2))^{p/2}\, \mu(ds)\Big)^{1/p}
= \tilde C_p \Big( \int_S (\E(|X(s)|^2))^{p/2}\, \mu(ds)\Big)^{1/p}\, .
\end{gather*}
the reverse inequality may be proved similarly. 

\medskip

The fact that a centered $X$ such that $\int_S (\E(|X(s)|^2))^{p/2}\, 
\mu(ds)<\infty$ 
is pregaussian follows from Lemma 5.1 of \cite{Rajput}. \hfill $\square$

\medskip

\subsection{Proof of Proposition \ref{necessity1}: the ASIP case.} 
Let $(X_n)_{n\ge 0}$ be i.i.d. variables in $L^2(\X)$. By assumption, 
they satisfy the ASIP. Hence, there exists i.i.d. gaussian variables 
$(W_n)_{n\ge 0}$, such that 
$$|X_0+\cdots +X_{n-1}-(W_0+\cdots +W_{n-1})|_\X=o(\sqrt{nL(L(n))})\qquad \as$$
Let $x^*\in \X^*$. By the law of the iterated logarithm 
(in the real case), $\E((x^*(X_0))^2)=\E((x^*(W_0))^2)$. In particular, 
$X_0$ is pregaussian. Then, we conclude thanks to Proposition 9.24 
of \cite{LTbook}. \hfill $\square$

\subsection{Proof of Proposition \ref{necessity2}: the BLIL case.}
Let $(X_n)_{n\ge 0}$ be i.i.d. pregaussian variables taking values in 
$\X$. By aasumption, they satisfy the BLIL. Let $1\le p<2$. It follows that 
$|X_n|_\X/n^{1/p}\underset{n\to +\infty}\longrightarrow 0$ $\P$-a.s. Hence, 
by the Borel-Cantelli lemma, $X_0\in L^p(\X)$. Then, the result follows 
from the proof of Proposition 9.25 of \cite{LTbook}.  
\hfill $\square$

\section{Proof of the martingale results}




\subsection{Proof of Proposition \ref{ineLIL}}
 This is just Proposition 3.3 of \cite{Cuny} 
 when $\X$ is $2$-smooth. Assume that $\X=L^p(S)$, $p\ge 1$. It suffices to prove the result 
when  $d\in L^p(S,L^2(\Omega,\F_0))$, otherwise $K_p(d)=+\infty$. There exists a sequence 
of step functions $(d_n)_{n\ge 1}$ converging in $L^p(S,L^2(\Omega,\F_0))$ 
to $d$. We may write $d_n(s,\omega)=\sum_{k=1}^{m_n} f_{k,n} (\omega) 
{\bf 1}_{A_{k,n}}(s)$, where $A_{k,n}\in \S$ and $f_{k,n} 
\in L^2(\Omega,\P)$. Let $\tilde d_n:= 
\sum_{k=1}^{m_n} (f_{k,n}-\E_{-1}(f_{k,n}))  {\bf 1}_{A_{k,n}}$. 
Then $(\tilde d_n)_{n\ge 1}$ converges to $d$ in $L^p(S,L^2(\Omega,\F_0))$
 as well 
(hence also in $L^2(\Omega,L^r(S))$, by Lemma \ref{grafakos}) 
and for every $s\in S$, $\tilde d_n(s,\cdot)$ is a real-valued 
martingale difference in $L^2(\Omega,\F_0,\P)$. Hence, applying  
Proposition \ref{ineLIL} to the $(2,1)$-smooth Banach space $\R$, we obtain that 
there exists $C_p>0$ such that for every $s\in S$,
\begin{equation}\label{ol}
\|\M_2(\tilde d_n(s,\cdot))\|_{p,\infty} \le C_{p} \|\tilde d_n(s,\cdot)\|_2\, .
\end{equation}
Notice that $\M_2(\tilde d_n, L^p(S))\le \Big( \int_S(\M_2(\tilde d_n(s,\cdot),
\R))^p  \, d\mu(s) 
\Big)^{1/p}$. Writing $\varphi(s,\cdot)=\M_2(\tilde d_n(s,\cdot),\R)$, it 
follows from lemma \ref{grafakos} that
$$
\| \M_2(\tilde  d_n, L^p(S))\|_{2,\infty}\le C_{p}
\Big( \int_S \|\varphi(s,\cdot)\|_{r,\infty}^p\, d\mu(s)\Big)^{1/p}\, .
$$
Then, we infer from \eqref{ol} that 
$$
\| \M_2(\tilde  d_n, L^p(S))\|_{r,\infty}\le C_{p}\Big(\int_S \|\tilde 
d_n(s)\|_2^p\, d\mu(s)\Big)^{1/p} \, .
$$
The desired result then follows by letting $n\to \infty$ (approximate 
first $\M_2$ by a supremum over a finite set of integers and use the 
monoton convergence theorem).  \hfill $\square$

\subsection{Proof of Proposition \ref{WIPmart}}

We shall first prove \eqref{inemart2} which will allow us to  derive the required  tightness  for the WIP. By Doob's maximal inequality for 
submartingales, we have 
\begin{gather*}
\|\max_{1\le k\le n} |S_k(d)|_\X\|_2 ^2
\le 2  \|\, |S_n(d)|_\X\|_2 ^2\, .
\end{gather*}
When $\X$ is $2$-smooth, \eqref{inemart2} then follows from 
\eqref{smooth} and the fact that, on Type 2 Banach spaces, the norms 
$\|\cdot \|_{\G(\X)}$ and $\|\cdot \|_{2,\X}$ are equivalent, by Proposition 9.24 of \cite{LTbook}. 

\smallskip

Assume now that $\X$ has cotype 2. Since $d$ is pregaussian, so is 
$S_n(d)$. Moreover, by orthogonality of real-valued martingale increments, 
we see that $G(S_n(d)/\sqrt n)=G(d)$. Since $\X$ has cotype 2, by Proposition 
9.25 of \cite{LTbook}, 
$$
\|S_n(d)\|_{2,\X}\le C \|G(S_n(d))\|_{2,\X}=C\sqrt n \|G(d)\|_{2,\X}
\le C\sqrt n \|d\|_{\G(\X)}\, ,
$$
and \eqref{inemart2} follows.

\smallskip

Let us prove the WIP. Let us recall 
the definition of tightness required here.

Let $X\in L^0(\Omega,\F_0,\P,\X)$. Recall that
$S_{n,t}=S_{n,t}(X):=S_{[nt]}+ (nt-[nt])X_{[nt]}$ and $T_{n,t}:=\frac{S_{n,t}
}{\sqrt n}$. 
We consider $\big( (T_{n,t})_{0\le t\le 1})_{n\ge 0}$ as a process taking 
values in $C([0,1],\X)$, the Banach space of continuous functions from 
$[0,1]$ to $\X$. 

\begin{defn}
We say that $\big( (T_{n,t})_{0\le t\le 1})_{n\ge 0}$ is tight if 
for every $\varepsilon>0$, there exists a compact set $\kappa$ 
of $C([0,1],\X)$ such that, 
$$
\P\big((T_{n,t})_{0\le t\le 1}\in \kappa\big) \ge 1-\varepsilon \qquad \forall 
n\ge 0.
$$
\end{defn}

Let $\X$ be either $2$-smooth or of cotype 2. Let $d\in \G(\X)$ 
with $\E_{-1}(d)=0$. Let us prove the  tightness of 
$\big((T_{n,t}(d))_{0\le t\le 1}\big)_{n\ge 1}$ 
in $C([0,1],\X)$. 

\medskip

We first recall the following   tightness criteria that may be easily deduced 
from Theorem 11.5.4 of Dudley \cite{Dudley}. 

\begin{lem}\label{tight}
Let $(\Gamma,\delta)$ be a separable complete metric space endowed with its Borel 
$\sigma$-algebra. Let $(\Omega,\F,\P)$ be a probability space and 
$(Z_n)_{n\ge 1}$ be a sequence of random variables on $\Omega$ taking values 
in $\Gamma$. Assume that, for every $\varepsilon >0$, there exist  $n_0\ge 1$
 and random variables $(Z_n^\varepsilon)_{n\ge n_0}$ such that  
\begin{itemize}
\item [$i)$] $(Z_n^\varepsilon)_{n\ge n_0}$ is tight; 
\item [$ii)$] $\sup_{n\ge n_0} \E(\delta(Z_n,Z_n^\varepsilon)) <\varepsilon$.
\end{itemize}
Then $(Z_n)_{n\ge 1}$ is tight.
\end{lem} 

Since $\X$ is separable, $\sigma(d)$ (the $\sigma$-algebra generated by $d$) 
is countably generated and there exists an increasing filtration 
$(\GG_m)_{m\ge 1}$ such that $\GG_m$ is finite for every $m\ge 1$ and 
$\sigma(d)=\vee_{m\ge 1} \GG_m$. For every $m\ge 1$, let $d_m:= \E(d|\GG_m)$. 
Since $\GG_m$ is finite, there exists $A_{1,m},\ldots , A_{k_m,m}\in \GG_m$ 
and $x_{1,m},\ldots , x_{k_m,m}\in \X$ such that $d_m=\sum_{1\le k\le k_m}
x_k {\bf 1}_{A_{k,m}}$. By Lemma \ref{martpreg}, $(d_m)_{m\ge 1}$ 
converges in $\G(\X)$ to $d$. Hence, writing $\tilde d_m:= 
d_m-\E_{-1}(d_m)$ and using Lemma \ref{LTcot}, $(\tilde d_m)_{m\ge 1}$ 
converges in $\G(\X)$ to $d$. 

\medskip

By the WIP for real-valued martingales with stationary and ergodic 
increments, for every $m\ge 1$, $((T_{n,t}(\tilde d_m))_{0\le 
t\le 1})_{n\ge 0}$ 
is tight in $C([0,1],\X)$. 

\smallskip

Now, by \eqref{inemart2},
\begin{gather*}
\|\sup_{0\le t\le 1}|T_{n,t}(\tilde d_m)-T_{n,t}(d)|_\X\|_2
\le \frac{3}{\sqrt n} \|\max_{1\le k\le n} |S_k(\tilde d_m)-S_k(d)|_\X\|_2 
\le C \|\tilde d_m-d\|_{\G(\X)} \underset{m\to \infty}\longrightarrow 0\, ,
\end{gather*}
and the tightness of $((T_{n,t}(d))_{0\le t\le 1})_{n\ge 0}$ in 
$C([0,1],\X)$ follows from Lemma \ref{tight}.

\medskip

Let us write $T_{n,t}(d)=T_{n,t}$. The second step consists in proving the convergence of the 
finite-dimensional laws. That is, it remains to prove that, 
for any $0= t_0< \ldots < t_m= 1$, $((T_{n,t_i}-T_{n,t_{i-1}})_{1\le i\le m})_{n\ge 1}$ 
converges in law to $(W_{t_i}-W_{t_{i-1}})_{1\le i\le m}$, where $(W_t)_{0\le t\le 1}$ 
is a brownian motion with covariance operator $\K_d$. Using tightness again (and the Cramer-Wold device), 
it suffices to prove that for any $0= t_0< \ldots < t_m= 1$ and any 
$x_1^*,\ldots , x_m^*\in \X^*$, $\sum_{i=1}^mx^*_i(T_{n,t_i}-T_{n,t_{i-1}})$ 
converges in law to $\sum_{i=1}^m x^*_i(W_{t_i}-W_{t_{i-1}})$ 
as $n\to \infty$. 

\medskip

Hence, we are back to prove a CLT for an array of martingale differences. 
Let us recall the following CLT of McLeish, as stated in Theorem 3.2 page 
58 of Hall and Heyde \cite{HH}. 

\begin{prop}
Let $(X_{n,j})_{1\le j\le k_n}$ be (real valued) martingale differences 
for every $n\ge 1$. Assume that there exists $\sigma\ge 0$ such that 
\begin{itemize}
\item [$(i)$] $\max_{1\le j\le k_n} |X_{n,j}|\overset{\P}\longrightarrow 0$;
\item [$(ii)$] $\sum_{1\le j\le k_n} X_{n,j}^2\overset{\P}\longrightarrow 
\sigma^2$; 
\item [$(iii)$] $\sup_{n\ge 1} \E(\max_{1\le j\le k_n}X_{n,j}^2)<\infty $. 
\end{itemize}
Then $(\sum_{1\le j\le k_n} X_{n,j})_{n\ge1}$ converges in law  to a normal 
law $\N(0,\sigma^2)$.  
\end{prop}

Take $k_n:=n$ and for every $1\le i\le m$ and every $[nt_{i-1}]\le j \le 
[nt_i]-1$, take $X_{n,j}:=x_i^*(d) \circ \theta^j/\sqrt n$. 

\medskip

Then, setting $Z:=\max_{1\le i\le m}|x_i^*(d)|$ (which belongs to $L^2(\Omega)$), 
we have   $\max_{1\le j\le k_n} |X_{n,j}|\le \max_{1\le j\le n} 
Z\circ \theta^j/\sqrt n$ which implies  $i)$, by the Borel-Cantelli lemma, and $iii)$ by standard arguments. Now, by the ergodic theorem we have 
$$
\frac{1}{n}\sum_{j=[nt_{i-1}]}^{[nt_i]} (x_i^*(d))^2\circ \theta^j
\underset{n\to \infty}\longrightarrow (t_i-t_{i-1}) \E((x_i^*(d)^2)\qquad \as\, ,
$$
 hence in probability. Hence the proof is complete. \hfill $\square$
 
 \medskip
 
 \subsection{Proof of Proposition \ref{ASIPmart}}
 
 Let us prove the CLIL. Notice that $\G_0(\X):= 
 \{d\in \G(\X,\F_0)~:~\E_{-1}(d)=0\}$ is a closed subspace of $\G(\X)$. 
 By \eqref{mm} and  Proposition \ref{banprin}, the set of $d\in \G_0(\X)$, such that 
 $(d\circ \theta^n)_{n\ge 0}$ satisfies the CLIL is closed in $\G_0(\X)$. 
 Then, the CLIL  follows by approximating any $d\in \G_0(\X)$ by 
 a martingale difference with values in a finite dimensional Banach space 
 as in the proof of Proposition \ref{WIPmart}.
 
 Then, \eqref{limsup} follows from a result of Kuelbs (see e.g. Proposition D. 
 of \cite{Cuny}) combined with the LIL for real valued 
 stationary (and ergodic) martingale differences. 
 
 \medskip
   To prove the ASIP, we just apply the following version of 
 Theorem 3.2 of Berger \cite{Berger} whose proof may be done 
 similarly. 
 
 \begin{theo}\label{Berger}
 Let $\X$ be a real separable Banach space. Assume that $\theta$ is ergodic. 
 Let $X\in L^0(\Omega,\F_0,\P,\X)$ 
 be such that   $\E(x^*(X)^2)<\infty$, for every $x^*\in \X^*$. 
 Assume that $(X\circ \theta^n)_{n\ge 0}$ satisfies the 
 CLIL and that for every $x^*\in \X^*$, there exists $Z=Z_{x^*}\in L^2(\Omega,\F_0,\R)$ 
 with $\E_{-1}(Z)=0$ such that  
 \begin{gather}\label{weakm}
 S_n(x^*(X))-S_n(Z) =o(\sqrt{nL(L(n))}) \qquad \as \\
 \|S_n(x^*(X))-S_n(Z)\|_2 =o(\sqrt n)\, .
 \end{gather}
Then, for every $x^*,y^*\in \X^*$, 
 $\K(x^*,y^*):=\lim_{n\to \infty} 
 \frac{{\rm cov}(x^*(S_n(X))y^*(S_n(X))}{n} $ exists. 
 Assume moreover that $\K$ is the covariance operator of a gaussian 
 variable. \\
 Then, $(X \circ \theta^n)_{n\ge 0}$ satisfies the ASIP.
 \end{theo}
 
 \section{Proof of the maximal inequalities}
 
 \subsection{Proof of Proposition \ref{inemax}}

We make the proof by induction. For $d=0$ we have 
$$S_1=X-\E_{-1}(X ) + \E_{-1}(X)= 
(X -\E_{-1}(X)) + \E_{-1}(S_1)$$
 and the result follows in that case. 
 
 \medskip
 
 Assume that we already proved the result for some $d\ge 0$. For every $1\le i\le 2^{d+1}$, we have
 \begin{gather*}
 S_i =\sum_{\ell=0}^{i-1} (X-\E_{-1}(X))\circ \theta^{\ell} + 
 \sum_{\ell=0}^{i-1} (\E_{-1}(X))\circ \theta^{\ell}\, ,
 \end{gather*}
and for every $1\le j\le 2^d$ (with $\sum_{\ell=0}^{-1}=0$),
\begin{gather*}
\sum_{\ell=0}^{2j-1} (\E_{-1}(X))\circ \theta^{\ell} 
=\sum_{\ell=0}^{j-1} (\E_{-1}(X)+ \E_{-1}(X)\circ \theta)
\circ \theta^{2\ell}\, ;\\
\sum_{\ell=0}^{2j-2} (\E_{-1}(X))\circ \theta^{\ell} 
=(\E_{-1}(X))\circ \theta^{2j-2}+ \sum_{\ell=0}^{j-2} (\E_{-1}(X)+ \E_{-1}(X)\circ \theta)\circ \theta^{2\ell}\, .
\end{gather*}

Hence,
\begin{gather}
\nonumber \max_{1\le i\le 2^{d+1}} |S_i|_\X \le \max_{1\le i\le 2^{d+1}}\Big|\sum_{\ell=0}^{i-1} (X-\E_{-1}(X))\circ \theta^{\ell } 
\Big|_\X + \max_{1\le j\le 2^d}|\E_{-1}(X)|_\X\circ \theta^{2j-2} \\
\label{ine0}+ 
\max_{1\le j\le 2^d}\Big| \sum_{\ell=1}^{j} (\E_{-1}(X)+ 
\E_{-1}(X)\circ \theta)\circ \theta^{2\ell}\Big|_\X \, .
\end{gather}

We shall apply our induction hypothesis to the following situation: 
$\tilde X:= \E_{-1}(X)+ 
\E_{-1}(X)\circ \theta$, the transformation 
$\tilde \theta:= \theta^2$ and the filtration given by $\tilde \F_n:=\tilde 
\theta^{-n}(\F)=\F_{2n}$ for every $n\in \Z$. 
\smallskip 

We shall also use the notation $\tilde \E_n(\cdot):= \E(\cdot |\tilde \F_n)$ 
and $\tilde S_n = \sum_{\ell=0}^{n-1} \tilde X\circ \tilde \theta^k$.

\smallskip

Notice then that we have 

\begin{gather*}
\tilde S_n  =\sum_{\ell=0}^{n-1} (\E_{-1}(X)+ 
\E_{-1}(X)\circ \theta)\circ \theta^{2\ell} , \qquad  
\tilde \E_{-2^k}(\tilde S_{2^k})= \E_{-2^{k+1}}(S_{2^{k+1}})\\
\quad \mbox{ and } \qquad \tilde X -\tilde \E_{-1}(\tilde X) 
=\E_{-1}(S_1)+\E_{-1}(S_1)\circ \theta -\E_{-2}(S_2)\, .
\end{gather*}

Hence, by our induction hypothesis and using the change of index $k\to k+1$, 
we infer that 

\begin{gather}\label{ine1}
\max_{1\le i\le 2^d} |\tilde S_i|_\X \le |\E_{-2^{d+1}}(S_{2^{d+1}})|_\X + \sum_{k=1}^{(d+1)-1}
\max_{0\le \ell \le 2^{(d+1)-1-k}-1}|\E_{-2^k}(S_{2^k})|_\X\circ \theta^{2^{k+1}\ell}\\
\nonumber 
+\sum_{k=1}^{(d+1)-1}\max_{1\le i\le 2^{(d+1)-k-1}}\Big|\sum_{\ell=0}^{i-1} 
\Big[ \E_{-2^k}(S_{2^k})+\E_{-2^k}(S_{2^k})\circ \theta^{2^k} 
-\E_{-2^{k+1}}(S_{2^{k+1}})\Big]\circ \theta^{2^{k+1}\ell}\Big|_\X \, .
\end{gather}

Then, the result follows by combining \eqref{ine0} and \eqref{ine1}. 
\hfill $\square$

\subsection{Proof of Corollary \ref{cormax2}} 
 We shall use Proposition 
\ref{inemax}. We first notice that 
$$
\max_{0\le \ell \le 2^{d-1-k}-1}|\E_{-2^k}(S_{2^k})|_\X\circ \theta^{2^{k+1}\ell}
\le \Big( \sum_{0\le \ell \le 2^{d-1-k}-1}|\E_{-2^k}(S_{2^k})|_\X^2\circ \theta^{2^{k+1}\ell} \Big)^{1/2}\, .
$$
Hence, using that $\theta$ preserves $\P$, we infer that
$$
\Big\| \max_{0\le \ell \le 2^{d-1-k}-1}|\E_{-2^k}(S_{2^k})|_\X\circ \theta^{2^{k+1}\ell} \Big\|_2\le 2^{(d-1-k)/2} 
\| \E_{-2^k}(S_{2^k})\|_{2,\X}\, .
$$

\medskip

Applying \eqref{inemart2} to (the martingale difference)  $d=X-\E_{-1}(X)$  we see that 
$$
\Big\|\max_{1\le i\le 2^d} \Big|\sum_{\ell=0}^{i-1} (X-\E_{-1}(X))\circ \theta^{\ell } \Big|_\X \Big\|_{\G(\X)} \le   C(\|X\|_{\G(\X)} 
+\|\E_{-1}(X)\|_{\G(\X)})\, .
$$

Similarly, we may apply \eqref{inemart2} with 
$d_k= \E_{-2^k}(S_{2^k})+\E_{-2^k}(S_{2^k})\circ \theta^{2^k} 
-\E_{-2^{k+1}}(S_{2^{k+1}})$ (and $\theta^{2^{k+1}}$ instead of $\theta$).
To conclude we just notice that, by Lemma 
\ref{LTcot}, $\|X-\E_{-1}(X) \|_{\G(\X)}\le (1+\sqrt 2)\|X\|_{\G(\X)}$ 
and that $\|\E_{-2^k}(S_{2^k})+\E_{-2^k}(S_{2^k})\circ \theta^{2^k} 
-\E_{-2^{k+1}}(S_{2^{k+1}})\Big]\circ \theta^{2^{k+1}\ell}\|_{\G(\X)} 
\le (1+\sqrt 2)^2\|\E_{-2^k}(S_{2^k})\|_{\G(\X)}$.  \hfill $\square$

\subsection{Proof of Proposition \ref{inemaxL2s} }

By Hopf's maximal inequality, for every $X\in L^1(\Omega,\R)$, and every 
measure preservint $\theta$
$$
\|\M_1(X, \theta)\|_{1,\infty}\le \|X\|_1\, .
$$

Then, the proposition follows from \eqref{inemax2} combined \eqref{mm}. 
\hfill $\square$

\section{Proof of the limit theorems under projective conditions}

Before doing the proof, let us give general facts about 
$\|\cdot \|_{MW_2}$, that will be used in the sequel.

\medskip 

 Define $MW_2:= \{X\in L^2(\Omega, \F_0,\P,\X)~:~\|X\|_{MW_2}<\infty\}$. Then, 
$(MW_2,\|\cdot \|_{MW_2})$ is a Banach space. 

\medskip

For every $X\in L^1(\Omega ,\F_0,\P,\X)$ define $QX =\E_0(X\circ \theta)$.  
Notice that $Q^n(X)=\E_0(X\circ \theta^n)$. Then, 
clearly $Q$ is a contraction of $L^2(\Omega,\F_0,\X)$ and, by Lemma 
\ref{LTcot}, $Q$ is power bounded on $\G(\X)$, i.e., for every 
$X\in \G(\X)$, $\sup_{n\ge 1}\|Q^nX\|_{\G(\X)}\le C \|X\|_{\G(\X)}$, for some 
universal $C>0$.
\medskip
Now, we see that 
\begin{gather*}
\|X\|_{MW_2}=\sum_{n\ge 0} \frac{\|\sum_{k=0}^{2^n-1}Q^k X\|_{\G(\X)}}
{2^{n/2}}\, .
\end{gather*}
Hence,  $Q$ is power bounded  on $MW_2$. 

\medskip

Writing $V_n:=I+\cdots +Q^{n-1}$ and using that $\|V_nV_kX\|_{\G(\X)}\le C\min(k\|V_n\|_{2,\X},n\|V_kX\|_{\G(\X)})$, we 
see that, for every $X\in MW_2$, 
\begin{equation}\label{mean}
\frac{\|V_{2^n} X\|_{MW_2}}{2^n } \le C\Big( \frac{\|V_{2^n}\|_{\G(\X)}}{2^{n/2}} 
+ \sum_{k\ge n+1} \frac{\|V_{2^k} X\|_{\G(\X)}}{2^{k/2}
}\Big)\underset{n\to+\infty} 
\longrightarrow 0\, .
\end{equation}
In particular, for every $m\ge 1$, taking $n$ such that $2^n\le m<2^{n+1}$, we have 
 $\|V_{m} X\|_{MW_2}\le C \sum_{k=0}^n \|V_{2^k}\|_{MW_2}=o(2^n)=o(m)$.

\medskip 

In particular, we see that $Q$ is mean ergodic on $MW_2$ and has no 
non trivial fixed point (see e.g. Theorem 1.3 p. 73 of \cite{Krengel}), i.e.,  
\begin{equation}\label{meanerg}
MW_2=\overline{(I-Q)MW_2}^{MW_2} \, .
\end{equation}

\subsection{Proof of Proposition \ref{propWIPcot} and 
Theorem \ref{theoASIPhilb}} In both results, $\X$ is a Banach space 
of cotype 2. Let $X\in (I-Q)MW_2$. Let $Y\in MW_2$ be the unique (notice 
that $Q$ has no fixed point on $MW_2$) solution 
to $X=(I-Q)Y$. Then, one may define 
$$\D(X):=Y-\E_{-1}(Y)=Y-QY\circ 
\theta^{-1}\, .
$$ 
Notice that $X=\D(X)+QY-QY\circ \theta^{-1}$ and that 
$\D(X)$ is a martingale difference. In particular 
\begin{equation}\label{marttt}
\|G(S_n(\D(X)))\|_{2,\X}
=\sqrt n \|G(\D(X))\|_{2,\X}\, .
\end{equation}

\medskip

Recall that, since $\X$ has cotype 2, there exists $C>0$, such that 
for every $Z\in \G(\X)$, 
\begin{equation}\label{G(X)}
\|G(Z)\|_{2,\X}/C\le \|Z \|_{\G(\X)}
\le C\|G(Z)\|_{2,\X} \, .
\end{equation}   

\medskip

Now, it follows from the proof of Proposition \ref{inemax} (combined with \eqref{marttt} applied to the martingales with stationary increments 
that appear in the proof) that 
there exists $D>0$ such that for every $d\ge 0$,
\begin{equation}\label{ed}
\|G(S_{2^d}(X))\|_{2,\X} \le 
D2^{d/2}\Big(  \|G(X)\|_{2,\X} + \sum_{k=0}^d 
2^{-k}\|G(\E_0(S_{2^k}(X))\|_{2,\X}\Big) \, .
\end{equation}

Notice that $\|S_{2^d}(QY-QY\circ \theta^{-1})\|_{\G(\X)} 
\le \|QY\circ \theta^{-1}\|_{\G(\X)}+\|QY\circ \theta^{2^d-1}\|_{\G(\X)}
=o(2^{d/2})$ and that $\|G\big( S_{2^d}(\D(X))\big)\|_{2,\X}\le 
\|G\big( S_{2^d}(X)\big)\|_{2,\X}+\|G(S_{2^d}(QY-QY\circ \theta^{-1})\|_{2,\X}$. 

Combining this with \eqref{ed}, \eqref{G(X)} and \eqref{marttt} and 
letting $d\to \infty$, we infer that 
$$
\|\D(X)\|_{\G(\X)}\le C \|X\|_{MW_2}\, .
$$
Hence, we may extend our linear operator $\D$ continuously to 
$\overline{(I-Q)MW_2}^{MW_2}=MW_2$. Notice that $\D$ takes values 
in $\G_0(\X)=\{Z\in \G(\X,\F_0)~:~\E_{-1}(Z)=0\}$.

\medskip

\smallskip

Let us prove Proposition \ref{propWIPcot}. By Corollary \ref{cormax2} 
and \eqref{inemart2}, there exists $C>0$ such that
\begin{equation*}\label{inemartMW}
\|\max_{1\le k\le n} |S_k(X)-S_k(\D(X))|_\X\|_2\le C \sqrt n\|X\|_{MW_2}\, .
\end{equation*}
By linearity of $\D$ (and of $X\mapsto S_k(X)$) it then suffices 
to prove \eqref{wipmax} for a set of $X$'s that is dense in 
$MW_2$, in particular for $X\in (I-Q)MW_2$. But if $X=(I-Q)Y$ with 
$Y\in MW_2$, we have, for every $K>0$ 
\begin{gather*}
\|\max_{1\le k\le n} |S_k(X)-S_k(\D(X))|_\X\|_2 \le 
\|\max_{1\le k\le n} |S_k(QY-QY\circ \theta^{-1})|_\X\|_2\\
\le \|QY\|_{2,\X}+ \|\max_{1\le k\le n} |QY\circ \theta^{k-1})|_\X\|_2
\le \|QY\|_{2,\X}+  K + n\| |QY|_\X {\bf 1}_{\{|QY|_\X\ge K\}}\|_2\, .
\end{gather*}
Hence 
\begin{gather*}
\limsup_{n\to \infty} \|\max_{1\le k\le n} |S_k(X)-S_k(\D(X))|_\X\|_2 
\le \| |QY|_\X {\bf 1}_{\{|QY|_\X\ge K\}}
\|_2\underset{K\to\infty}\longrightarrow 0 \, ,
\end{gather*} 
and \eqref{wipmax} holds. Then, the proof of the WIP follows 
 from Lemma \ref{tight} and Proposition \ref{WIPmart}.

\medskip

Let us prove Theorem \ref{theoASIPhilb}. By Proposition \ref{ineLIL} and 
\eqref{inemaxL2}, for every $1<p<2$, there exists $C_p>0$ such that 
$$
\|\M_2(X-\D(X))\|_{p,\infty}\le C_p \|X\|_{MW_2}\, .
$$
Hence, by the Banach principle, see Lemma \ref{banprin}, it suffices 
to prove \eqref{asipmax} for $X=(I-Q)Y$, with $Y\in MW_2$. 
But in this case the result is obvious, since $|QY|_\X\in L^2(\Omega)$ 
and, by the Borel-Cantelli lemma, $|QY|_\X\circ \theta^{n-1}=o(\sqrt n)$ \as 
\, By \eqref{asipmax} and Proposition \ref{ASIPmart}, 
$(X\circ \theta^n)_{n\ge 0}$ satisfies the CLIL. Then, the ASIP follows from Proposition \ref{Berger}, using that $\D(X)$ is pregaussian. 

\medskip

It remains to prove \eqref{kuelbslimsup}. The first equality follows from 
\eqref{asipmax} and \eqref{limsup}. Let us prove that, with 
$d=\D(X)$,  $\sup_{x^*\in \X^*, 
|x^*|_{\X^*}\le 1}\|x^*(d)\|_2\le 10\sqrt 2 \|X\|_{MW_2}$.  We first notice that 
$x^*(d)=\D(x^*(X))$ (with the obvious "new" meaning of the operator $\D$).  
Proceeding as above one can prove that for every $m\ge 0$, 
$$
\|x^*(d)\|_2= 2^{m/2} \|S_{2^m}(d)\|_2/2^{m/2}\le \|S_{2^m}(X)\|_2/2^{m/2}+ 
\|S_{2^m}(d)-S_{2^m}(X)\|_2/2^{m/2}\, .
$$ 
Applying Proposition \ref{propWIPcot} (noticing that 
$\|x^*(X)\|_{MW_2}\le \|X\|_{MW_2}$) and Corollary \ref{cormax2} to $x^*(X)$, 
  we derive that $\|x^*(d)\|_{MW_2}\le 10 \sqrt 2\|X\|_{MW_2}$ 
  and the proof is complete.

\hfill $\square$

\subsection{Proof of Theorem \ref{theoASIPsmooth}} 

Let us prove the WIP. As above we shall first prove tightness. 
Let $X\in MW_2$. Let $\varepsilon >0$. By \eqref{meanerg}, there 
exists $Y\in MW_2$ such that $\|X-(I-Q)Y\|_{MW_2}\le \varepsilon$. 

\smallskip

Then, by Corollary \ref{cormax2}, 
\begin{gather*}
\|\max_{1\le k\le n} |S_k(X)-S_k((I-Q)Y)|_\X\|_2 \le 
 C \varepsilon \sqrt n.
\end{gather*}

Now, as in the proof of Proposition \ref{propWIPcot}, for every 
$K>0$ we have
\begin{gather*}
\|\max_{1\le k\le n} |S_k((I-Q)Y)-S_k(Y-\E_{-1}(Y)))|_\X\|_2 \le 
 \|QY\|_{2,\X}+  K + n\| |QY|_\X {\bf 1}_{\{|QY|_\X\ge K\}}\|_2\, .
\end{gather*}
Chose  $K$  such that $\| |QY|_\X {\bf 1}_{\{|QY|_\X\ge K\}}\|_2\le 
\varepsilon $ and then chose  $n_0\ge (\|QY\|_{2,\X}+  K)^2/\varepsilon ^2$.

\smallskip
Then, $\|\sup_{0\le t\le 1} |T_{n,t}(X)-T_{n,t}(Y-\E_{-1}(Y))|_\X
\le C \varepsilon$. Now, $Y-\E_{-1}(Y)$ is a martingale difference, hence, 
by Proposition \ref{WIPmart}, $\big(( T_{n,t}(Y-\E_{-1}(Y))_{0\le t\le 1}\big){n\ge 0}$ is tight in $C([0,1],\X)$. Then, the tightness of 
$\big(( T_{n,t}(X)_{0\le t\le 1}\big)_{n\ge 0}$ follows from Proposition 
\ref{tight}. 

\medskip

The proof of the finite-dimensional laws may be done exactly 
as the proof of the martingale case, hence is ommitted. 
The fact that the covariance operator is given as stated follows 
from the fact that for any $x^*\in \X^*$, 
$x^*(X)$ satisfies the assumption of Proposition \ref{propWIPcot}.

\medskip

Let us prove the ASIP. We shall use Theorem \ref{Berger}. In 
particular, we have to prove that $(X\circ \theta^n )_{n\ge 0}$ 
satisfies the CLIL. 

\smallskip

By \eqref{inemaxL2} and Lemma \ref{banprin}, the set 
$\{X\in MW_2~:~(X\circ \theta^n)_{n\ge 0} \mbox{ statisfies the CLIL}\}$ 
is closed in $MW_2$. Hence, it suffices to prove the 
CLIL for $X=(I-Q)Y$, with $Y\in MW_2$. But then, $X=Y-\E_{-1}(Y)+ QY\circ \theta^{-1} 
-QY$ and $((Y-QY)\circ \theta^n)_{n\ge 0}$ satisfies the CLIL by 
Proposition \ref{ASIPmart}, while $|S_n( QY\circ \theta^{-1} -QY)|_\X
=o(\sqrt n)$ \as, by the Borel-Cantelli lemma. Hence the CLIL is proved. 

\medskip

 \medskip
 
 Now, let $x^*\in \X^*$. Clearly, $x^*(X)$ satisfies the assumption of 
 Theorem \ref{theoASIPhilb}, taking for $\X$ the Hilbert space $\R$. 
 In particular, there exists $Z\in L^2(\Omega,\F_0,\R)$ with $\E_{-1}(Z)=0$ 
 such that
  \begin{gather}\label{weakm}
 S_n(x^*(X))-S_n(Z) =o(\sqrt{nL(L(n))}) \qquad \as \\
 \|S_n(x^*(X))-S_n(Z)\|_2 =o(\sqrt{n})\, .
 \end{gather}
  The fact that  
 $\K(x^*,y^*):=\lim_{n\to \infty} 
 \frac{{\rm cov}(x^*(S_n(X))y^*(S_n(X))}{n} $ is the covariance operator of a gaussian variable, follows from the WIP. 
 
 \medskip
 
 To prove the equality in \eqref{kuelbslimsup2}, by a result of Kuelbs (see e.g. Proposition 
 D.1 in \cite{Cuny}), we have to prove that for every 
 $x^*\in \X^*$, we have 
 $$
\limsup_{n} \frac{S_n(x^*(X))}{\sqrt{2nL(L(n))} }= \big(\K(x^*,x^*)\big)^{1/2} 
\qquad \as 
 $$
 But this follows from Theorem \ref{theoASIPhilb} applied to $x^*(X)$. Then, the inequality 
 in \eqref{kuelbslimsup2} may be proved as the inequality in 
 \eqref{kuelbslimsup}.  \hfill $\square$
 
\subsection{Proof of Proposition \ref{prop-PU}} We first recall tha construction of Peligrad and Utev \cite{PU}. 

We consider the  Markov chain 
$(W_n)_{n\ge 0}$ with state space ${\mathbb N}:=\{0,1,\ldots\}$ and transition 
probability given by $p_{i,i-1}=1$ and  $p_{0,i-1}=p_i
=\P(\tau=i)$ for every $i\ge 1$, and 
$p_{i,j}=0$ otherwise. The stationarity is guaranteed by the 
condition $\E(\tau)<\infty$ and then, the stationary distribution 
  $\pi:=(\pi_i)_{i\ge 0}$ is given by $\pi_0=1/\E(\tau)$ and 
$\pi_i=\pi_0\sum_{j\ge i+1}p_j$. 

\medskip

Since our Markov chain is stationary, we may consider its two-sided version 
$(W_n)_{n\in \Z}$, taking for $(\Omega,\F,\P)$ the canonical space, for 
$\theta$ the shift and for $\F_0$, $\sigma\{W_n~:~n\le 0\}$. Then we are exactly in the situation considered in our paper.

\medskip

Let $(a_n)_{n\in \N}$ be a sequence of positive numbers with 
$a_n\to 0$ as $n\to \infty$. It is proved in \cite{PU} that there exists a choice of $(p_n)_{n\ge 0}$, such that $\E(\tau)<\infty$, 
$\E(\tau^2)=+\infty$ and such that \eqref{MWPU} holds with 
$X:={\bf 1}_{\{W_0=0\}}-\pi_0$.

Define $b_n:=\sqrt{n\log \log n}$. Let us prove that $\limsup_n |S_n|/b_n=+\infty$ \as  

\smallskip

Let $T_0:=0$ and, for $k\ge 1$, $T_k:=\min\{t>T_{k-1}~:~W_t=0\}$. 
Define then, $\tau_k:=T_k-T_{k-1}$. Then, $(\tau_k)_{k\ge 1}$ is iid,
distributed like $\tau$ and $S_{T_k}=\sum_{i=1}^k(1-\pi_0\tau_i)$.

\medskip 

It is enough to prove that $\limsup_k |S_{T_k}|/b_{T_k} 
=+\infty$ \as

\medskip

Since $\E(\tau)<\infty$, by the strong law of large numbers, 
$T_n/n\underset{n\to \infty}\longrightarrow \E(\tau)$ \as, hence it is enough 
to prove that $\limsup_k |S_{T_k}|/b_{k} 
=+\infty$ \as \, In particular, it is enough to prove that 
\begin{equation}\label{ii}
\limsup_k \big|\sum_{i=0}^k (1-\pi_0\tau_i)\big|/b_{k} 
=+\infty\qquad \as
\end{equation} 

But \eqref{ii} follows from Strassen's converse to the law of the iterated logarithm, see for instance \cite{LTbook} page 203-204, since 
$\E(\tau^2)=+\infty$. \hfill $\square$

\section{Technical results}

We recall here the Banach principle that we need (see Proposition C.1 
of \cite{Cuny}).

\begin{lem}\label{banprin}
Let $(\Omega,\F,\P)$  be a probability space and $\X, \B$ be
Banach spaces. Let $\C$ be a vector space of measurable
functions from $\Omega$ to $\X$. Let $(T_n)_{n\ge 1}$ be a sequence
of linear maps from $\B$ to $\C$. Assume that there exists
a positive decreasing function $L$ on $]0,+\infty[$, with
$\lim_{\lambda \to \infty}L(\lambda)=0$, such that
\begin{equation}\label{bancond}
\P(\sup_{n\ge 1} |T_n x|_\X> \lambda |x|_\B)\le L(\lambda) \qquad
\forall \lambda >0,x\in \B\, .
\end{equation}
Then the set $\{x\in \B \,:\,(T_n x)_{n\ge 1} \mbox{ is \as \,relatively compact
in $\X$ }\}$
and the set $\{x\in \B\,:\, |T_n x|_\X\to 0 \quad \as\}$
are closed in $\B$.
\end{lem}

We give here a technical result concerning $L^r$ spaces of $L^p$-valued 
variables.

\begin{lem}\label{grafakos}
Let $1\le p < r<\infty$. Let $(\Omega,\F,\P)$  be a probability space 
and $(S,\S,\nu)$ 
be a $\sigma$-finite measure spaces. There is a continuous embedding 
from $L^p(S,L^{r,\infty}(\Omega))$ (resp. $L^p(S,L^{r}(\Omega))$) into 
$L^{r,\infty}(\Omega,L^{p}(S))$ (resp. $L^{r}(\Omega,L^{p}(S))$).
\end{lem} 
\noindent {\bf Proof.} We first recall some useful fact about weak $L^r$-spaces (see Exercise 
1.1.11 p. 13 of Grafakos \cite{Grafakos}).  
For every $r>1$ and every $0<t<r$, let 
$$N_{r,t}(|X|_\X) := 
\sup_{\P(A)>0}\frac{1}{\P(A)^{1/r-1/t}}\Big( \E( |X|_\X^t{\bf 1}_A)\Big)^{1/t}\, .
$$ 
Then, there exists $C_{r,t}$ such that 
$$
\|X\|_{r,\infty,\X}/C_{r,t}\le N_{r,t}(|X|_\X)\le C_{r,t}\|X\|_{r,\infty,\X} \, ,
$$ 
and for $t=1$, $N_{s,1}$ is a norm.
\medskip 

Let $f(s,\omega)= \sum_{i=1}^n f_i(\omega)
{\bf 1}_{A_i(s)}$ 
be a step function of  $L^p(S,L^{r,\infty}(\Omega))$, i.e. $A_i\in \S$ and 
$f_i\in L^{r,\infty}(\Omega)$. We may consider $f$ as an element of 
$L^0(S\times\Omega, \S\otimes \F)$ or as an element of 
$L^0(\Omega,\F, L^0(S,\S))$.

\medskip 

Take $\X=L^p(\mu)$ and $t=p$. We have, using Fubini, 
\begin{gather*}
\E(\|f\|_{L^p(\mu)}^p{\bf 1}_A)= \int_S \E(|f(s,\cdot)|^p{\bf 1}_A)\, d\mu(s)
\le \P(A)^{1/r-1/p} \int_{S}N_{r,p}(|f|(s,\cdot))\, d\mu(s)\, .
\end{gather*}
Hence, 
$$
\|f\|_{r,\infty,L^p(S)}\le C_{r,p}^2 \Big( \int_S \|f(s,\cdot)\|_{r,\infty}^p
\, d\mu(s)\Big)^{1/p}\, .
$$
 Hence, the identity map  sends step functions of $L^p(S,L^{r,\infty}(\Omega))$ 
 to elements of $L^{r,\infty}(\Omega,L^{p}(S))$ in a continuous way. 
 In particular, it can be extended continuously in an injective map to the 
 whole  $L^p(S,L^{r,\infty}(\Omega))$.\hfill $\square$
 
 \end{appendix}

\end{document}